\documentclass[12pt]{article}

\usepackage{subcaption}

\usepackage{float}

\usepackage{multirow}

\usepackage{tikz}

\usetikzlibrary{calc,backgrounds,arrows,matrix}

\usepackage{enumerate}

\usepackage{color}
\definecolor{lightblue}{rgb}{0,0.2,0.5}

\usepackage{xcolor}
\definecolor{ForestGreen}{RGB}{34,139,34}
\definecolor{mauve}{rgb}{0.7,0,0.43}
\definecolor{dkgreen}{rgb}{0,0.6,0}
\definecolor{darkgreen}{rgb}{0,0.6,0}
\definecolor{darkorange}{rgb}{1.0, 0.55, 0.0}
\definecolor{lightblue}{rgb}{0,0.2,0.5}
\definecolor{blue1}{rgb}{0,0.1,0.9}

\usepackage[colorlinks=true, urlcolor=blue,linkcolor=blue, citecolor=lightblue]{hyperref}

\usepackage{listings}
\lstdefinelanguage{Maple}{
    morekeywords={proc, if, return, map, op, int, for, do, local, nops, convert, end},
    sensitive=false, 
    morecomment=[l]{//}, 
    morecomment=[s]{/*}{*/}, 
    morestring=[b]" 
} 

\usepackage{amssymb,amsmath}

\usepackage{graphicx}

\usepackage{graphics}
\usepackage{color}

\usepackage{tikz}
\usetikzlibrary{trees}

\DeclareMathAlphabet{\eufrak}{U}{}{}{}
\SetMathAlphabet\eufrak{normal}{U}{euf}{m}{n}
\SetMathAlphabet\eufrak{bold}{U}{euf}{b}{n}

\allowdisplaybreaks

 \def\qu{{\mathord{\mathbb Z}}}

 \def\Var{{\mathrm{{\rm Var}}}}

 \def\T{{\mathrm{{\rm T}}}}

 \def\sZZ{{\rm Z\kern-.45em{}Z}}

 \def\sQQ{{\kern 0.27em \vrule height1.45ex width0.03em depth0em
           \kern-0.30em \rm Q}}
 \def\qu{{\mathchoice
         {\sQQ}
         {\sQQ}
   {\kern 0.225em \vrule height1.05ex width0.025em depth0em \kern-0.25em \rm Q}
   {\kern 0.180em \vrule height0.78ex width0.020em depth0em \kern-0.20em \rm Q}
         }}
 \def\sGG{{\kern 0.27em \vrule height1.45ex width0.03em depth0em
           \kern-0.30em \rm G}}
 \def\gg{{\mathchoice
         {\sGG}
         {\sGG}
   {\kern 0.225em \vrule height1.05ex width0.025em depth0em \kern-0.25em \rm G}
   {\kern 0.180em \vrule height0.78ex width0.020em depth0em \kern-0.20em \rm G}
         }}

 \newtheorem{prop}{Proposition}[section]
 \newtheorem{lemma}[prop]{Lemma}
 \newtheorem{definition}[prop]{Definition}
 
 \newtheorem{theorem}[prop]{Theorem}
 \newtheorem{remark}[prop]{Remark}

\numberwithin{equation}{section}

\newcommand{\re}{\mathrm{e}}

 \newcounter{hyp}
\usepackage{aeguill}

\newenvironment{Proof}{\removelastskip\par\medskip \noindent{\em Proof.} \rm}{\penalty-20\null\hfill$\square$\par\medbreak}

\def\bprf{\begin{Proof}}
\def\nprf{\end{Proof}}
\def\bdes{\begin{description}}
\def\ndes{\end{description}}

\newtheorem{thm}{Theorem}[section]

\def\bdef{\begin{defn}}
\def\ndef{\end{defn}}
\def\bthm{\begin{thm}}
\def\nthm{\end{thm}}
\def\bprop{\begin{prop}}
\def\nprop{\end{prop}}
\def\brmk{\begin{remark}}
\def\nrmk{\end{remark}}
\def\bexa{\begin{exa}}
\def\nexa{\end{exa}}
\def\blem{\begin{lem}}
\def\nlem{\end{lem}}
\def\bcor{\begin{cor}}
\def\ncor{\end{cor}}
\def\bexe{\begin{exe}}
\def\nexe{\end{exe}}

\newcommand{\E}{\mathbb{E}}

\newcommand{\real}{\mathbb{R}}

\def\Var{\mathop{\hbox{\rm Var}}\nolimits}

\def\og{\leavevmode\raise.3ex
     \hbox{$\scriptscriptstyle\langle\!\langle$~}}
\def\fg{\leavevmode\raise.3ex
     \hbox{~$\!\scriptscriptstyle\,\rangle\!\rangle$}~}

\title{\Huge
 Numerical evaluation of ODE solutions by Monte Carlo enumeration of Butcher series 
}

\author{
  Guillaume Penent\footnote{\href{mailto:PENE0001@e.ntu.edu.sg}{pene0001@e.ntu.edu.sg}}
  \qquad
      Nicolas Privault\footnote{
\href{mailto:nprivault@ntu.edu.sg}{nprivault@ntu.edu.sg}
}
  \\
\small
Division of Mathematical Sciences
\\
\small
School of Physical and Mathematical Sciences
\\
\small
Nanyang Technological University
\\
\small
21 Nanyang Link, Singapore 637371
}

\allowdisplaybreaks

\usepackage{empheq}

\usepackage{bbm}

\makeatletter
\newcommand*\rel@kern[1]{\kern#1\dimexpr\macc@kerna}
\newcommand*\widebar[1]{
  \begingroup
  \def\mathaccent##1##2{
    \rel@kern{0.8}
    \overline{\rel@kern{-0.8}\macc@nucleus\rel@kern{0.2}}
    \rel@kern{-0.2}
  }
  \macc@depth\@ne
  \let\math@bgroup\@empty \let\math@egroup\macc@set@skewchar
  \mathsurround\z@ \frozen@everymath{\mathgroup\macc@group\relax}
  \macc@set@skewchar\relax
  \let\mathaccentV\macc@nested@a
  \macc@nested@a\relax111{#1}
  \endgroup
}
\makeatother

\usepackage{array}

\usepackage{relsize} 

\newcolumntype{C}[1]{>{\centering\let\newline\\\arraybackslash\hspace{0pt}}m{#1}}

\begin{document}

\maketitle

\baselineskip0.6cm

\vspace{-0.6cm}

\begin{abstract}
We present an algorithm for the numerical solution of ordinary differential equations by random enumeration of the Butcher trees used in the implementation of the Runge-Kutta method. Our Monte Carlo scheme allows for the direct numerical evaluation of an ODE solution at any given time within a certain interval, without iteration through multiple time steps. In particular, this approach does not involve a discretization step size, and it does not require the truncation of Taylor series. 
\end{abstract}



\baselineskip0.65cm

\section{Introduction}
Butcher series \cite{butcher1963},
\cite{butcherbk} are a powerful tool used to represent 
 the Taylor expansions appearing in the Runge-Kutta methods for the numerical
solution of ordinary differential 
equations (ODEs), see Chapters~4-6 of \cite{deuflhard}, and
 \cite{mclachlan} for a recent review 
 starting from the early work of \cite{cayley}. 
 Those series are making use rooted tree enumeration,
 which have many applications ranging from 
 geometric numerical integration to 
 stochastic differential equations, see for instance \cite{ehairer}
 and references therein, and 
 \cite{gubinelli}, \cite{bruned}, \cite{fossy} for the use of
 decorated trees for stochastic partial differential equations 
 and their connections with 
 the Butcher-Connes-Kreimer Hopf algebra
 \cite{connes}.

\medskip 

 It is known that the solution $y(t)$ of the autonomous
 $d$-dimensional ODE system 
\begin{equation}
\label{fjklds}
\begin{cases}
y'(t) = f(y(t))\\
y(0) = y_0\in \real^d, \qquad t \in \real_+, 
\end{cases}
\end{equation}
 where $f(y)=(f_1(y),\ldots , f_d(y))$ is a smooth 
 $\real^d$-valued function of $y$ in a domain of $\real^d$, can
 be expressed as
\begin{equation}
\label{sdfh} 
y(t) = y_0 + t f(y_0) + \frac{t^2}{2} f'(y_0)f(y_0) + \frac{t^3}{6} f'(y_0)f(y_0)'f(y_0)
+ \frac{t^3}{6} f''(y_0)[f(y_0),f(y_0)]
+ \cdots
\end{equation}
where we use the notation
\begin{align}
  \label{aa1}
  & 
f'(y_0)f(y_0) := \left(  \sum_{i_2=1}^d \frac{\partial f_{i_1}}{\partial x_{i_2}}(y_0) f_{i_2}(y_0)\right)_{i_1=1,\ldots , d},
\\
\nonumber
\\
\label{aa2}
  & f'(y_0)f'(y_0)f(y_0)
:= \left( \sum_{i_2,i_3=1}^d \frac{\partial f_{i_1}}{\partial x_{i_2}}(y_0)
\frac{\partial f_{i_2}}{\partial x_{i_3}}(y_0)
f_{i_3}(y_0)\right)_{i_1=1,\ldots , d},
\\
\nonumber
\\
\label{aa3}
  &  f''(y_0)[f(y_0),f(y_0)]
:= \left( \sum_{i_2,i_3=1}^d \frac{\partial^2 f_{i_1}}{\partial x_{i_2}\partial x_{i_3}}(y_0)
f_{i_2}(y_0)f_{i_3}(y_0)\right)_{i_1=1,\ldots , d},
\end{align} 
 etc.
 In addition, the expansion \eqref{sdfh} can be coded and enumerated using
 the following sequence of Butcher trees, which can also be
 represented without referring to the derivatives of $f$,
 see for example Table~1.1 page 53 of \cite{ehairer}. 
 
             \begin{table}[H]
  \centering
 \resizebox{1\textwidth}{!}
           {
  \begin{tabular}{||c|c|C{2.4cm}||c|c|C{2.4cm}||c|c|C{2.4cm}||}
 \hline
 Order & Coefficient & Butcher tree &  Order & Coefficient & Butcher tree & Order & Coefficient & Butcher tree  \\
 \hline \hline
 0 & $y_0$ & 
 $ \emptyset$ & 1 & $ f(y_0) $ &  
\scalebox{0.8}
       {
         \begin{tikzpicture}[grow=right, sloped]
\node[circle,draw,black,text=black,thick,minimum size=0.2cm]{$\small f$};
\end{tikzpicture}
}
 & 2      &
       $f'f(y_0)$ &
\scalebox{0.8}
       {
\begin{tikzpicture}[grow=up, sloped]
\node[circle,draw,black,text=black,thick]{$f'$}
    child[draw=blue]{
        node[circle,draw,blue,text=blue,thick] {$f$}
            edge from parent
    };
\end{tikzpicture}}
\\
\hline
\hline
 3 & $f'f' f(y_0)$ &
\scalebox{0.7}
    {
\begin{tikzpicture}[grow=up, sloped, node distance=0.1cm]
\node[circle,draw,black,text=black,thick]{$f'$}
    child[draw=blue] {
        node[circle,draw,blue,text=blue,thick] {$f'$}
        child[draw=purple] {
            node[circle,draw,purple,text=purple,thick] {$f$}
                edge from parent
        edge from parent
    }
    }
    ;
\end{tikzpicture}}
    &
     3 & $f''[f,f](y_0)$ &
\scalebox{0.8}
       {
\begin{tikzpicture}[grow=up, sloped]
\node[circle,draw,black,text=black,thick]{$f''$}
    child[draw=purple] {
        node[circle,draw,purple,text=black,thick] {$\textcolor{purple}{f}$}
            edge from parent
    }
    child[draw=blue] {
        node[circle,draw,blue,text=black,thick] {$\textcolor{blue}{f}$}
            edge from parent
    }
    ;
\end{tikzpicture}}
 \\
\cline{1-6}
\cline{1-6}
\end{tabular}
           }
\end{table}
                      
\vspace{-0.4cm}
 
\noindent
The numerical evaluation of Butcher series involves tree enumeration up to a certain
order that determines the level of precision of the algorithm.

\medskip 

In this paper, we propose an alternative method to the numerical
evaluation of ODE solutions,
 based on a random enumeration of Butcher trees by Monte Carlo simulation. 
 Probabilistic methods
 based on the Feynman-Kac formula
 provide alternatives to finite difference schemes,
 and have been successfully applied to the solution of
 partial differential equations. 
In particular, stochastic branching mechanisms
have been used to represent the solutions of partial differential
equations in \cite{skorohodbranching},
\cite{inw}, \cite{N-S}, \cite{hpmckean}, 
 \cite{lm}, \cite{chakraborty}.
 This branching argument has been recently extended in
 \cite{labordere} to the treatment of polynomial non-linearities in
 gradient terms, see also
 \cite{penent} 
 for nonlocal and fractional PDEs.

\medskip

In Theorem~\ref{t1}, under suitable integrability conditions
we express ODE solutions as the expected value of a functional of
 random Butcher trees which encode nonlinearities. 
 Then in Proposition~\ref{p1}
 we provide sufficient conditions
 ensuring that the representation formula of
 Theorem~\ref{t1} holds at any time within certain interval.
 Numerical values of ODE solutions can be computed beyond that
 initial interval by iterating the method and by piecing together
 the solutions obtained on adjacent intervals. 
 As noted before Proposition~\ref{p1},
 the integrability conditions are stronger in higher dimensions.   

\medskip

 This approach complements the use of the
 Feynman-Kac formula for the numerical estimation of the solutions
 of partial differential equations,
 see also \cite{selk} for a different approach to the
 Feynman-Kac representation of ODE solutions.
 Other links between Butcher trees and probability theory have been
 pointed out in \cite{mazza}, 
 see also \cite{skilling}
 for the numerical solution of ODEs as an inference problem by Bayesian techniques. 

\medskip

In comparison with integrators in the Runge-Kutta method,
our method requires the evaluation of $f$ and its
partial derivatives up to any order, whereas the
Runge-Kutta method uses only $f$.
The complexity of our algorithm grows linearly with dimension, as 
$d$ trees are used to solve a $d$-dimensional ODE system,
whereas the complexity of finite difference methods is generally polynomial
in the dimension $d$, depending on the order chosen in the truncation of
\eqref{sdfh}. 
Complexity in time can be estimated via the mean length of binary trees, 
which grows exponentially in time independently of dimension $d\geq 1$
in the case of exponentially distributed branch lifetimes,
as noted in Section~\ref{s3}. 
On the other hand, our method can be used to exactly approximate the solution
at any given time on a (possibly infinite) time interval. 
Our approach also benefits from the advantages of Monte Carlo estimators
whose computation can be paralleled straightforwardly.

 \medskip

 In Section~\ref{s2} we introduce the construction of coding trees that
 will be used for the numerical solution of ODEs. 
 Section~\ref{s3} presents the probabilistic representation formula
 of ODE solutions obtained by the random generation of coding trees. 
 In Section~\ref{s4} we consider examples and in Section~\ref{s5} we
 describe the correspondance between Butcher trees
 and coding trees, namely we show how any Butcher tree
 can be recovered by performing a depth first search
 on the corresponding coding tree, showing how Butcher series can be
 rewritten as series of expected values. 
 Section~\ref{s6} considers numerical applications, and the appendix
 contains the corresponding computer codes in Maple, Mathematica and Python. 
\section{Codes and mechanism} 
\label{s2}
 This section introduces the random coding trees used for the
 probabilistic representation of ODE solutions.
 We consider a multidimensional autonomous system of the form 
\begin{equation}
\label{Eautsys}
\left\{ 
\begin{array}{l}
  y_1(t) =  \displaystyle y_1(0) + \int_0^t f_1(y_1(s), \ldots, y_d(s)) ds\\
y_2(t) =  \displaystyle y_2(0) + \int_0^t f_2(y_1(s), \ldots, y_d(s)) ds\\
  \qquad \ \ \vdots \\
y_d(t) =  \displaystyle y_d(0) + \int_0^t f_d(y_1(s), \ldots, y_d(s)) ds, 
\end{array}
\right. 
\end{equation}
 $t\in \real_+$, where $f_i$ is a smooth Lipschitz function defined on a domain of $\real_+ \times \real^d$, $i=1,\ldots , d$.
 In order to solve \eqref{Eautsys} iteratively,
 we can start by expanding $f_i(y_1(s),\ldots , y_d(s))$ as 
\begin{equation}
\label{fjdkslf} 
f_i(y_1(s),\ldots , y_d(s)) = f_i(y_1(0),\ldots , y_d(0)) +
\sum_{j=1}^d \int_0^s f_j(y_1(u),\ldots , y_d(u)) \partial_j f_i(y_1(u),\ldots , y_d(u)) du 
\end{equation}
by differentiating $v(s) := f_i(y_1(s),\ldots , y_d(s))$, 
where we use the notation $\partial_j f_i = \partial f_i / \partial y_j$,
$i,j=1,\ldots , d$. 
 In the sequel, given $g$ a function from $\real^d$ into $\real$, 
 we let $g^*$ denote the mapping 
 \begin{align}
\nonumber 
    g^* : (\real^d)^{\real_+} & \longrightarrow \real^{\real_+}
    \\
    \label{fdhskjfd}
    (t\mapsto y(t)) & \longmapsto g^*(y):=(t \mapsto g(y(t))), 
\end{align} 
where $(\real^d)^{\real_+}$ represents the set of functions from $\real_+$ to $\real^d$.
In order to formalize and extend the iteration initiated in \eqref{fjdkslf},
we introduce the following definitions. 
In the sequel we let ${\rm Id}_i$ denotes the $i$-$th$ canonical projection
from $(\real^d)^{\real_+}$ to $\real^{\real_+}$ with 
 ${\rm Id}_i(y_1,\ldots,y_d) = y_i$, $i=1,\ldots , d$. 
\begin{definition}
 We let $\mathcal{C}$ denote the set of 
 functions from $(\real^d)^{\real_+}$ to $\real^{\real_+}$
 called \textit{codes}, defined as  
$$
\mathcal{C} := \big\{ {\rm Id}_i, \ \big( \partial_1^{i_1}\cdots \partial_d^{i_d} f_i \big)^* \ : \ i_1,\ldots , i_d \geq 0, \ i =1, \ldots, d \big\} 
$$
\end{definition}
 By \eqref{fdhskjfd}, the elements of $\mathcal{C}$ are operators
 mapping a function $h\in (\real^d)^{\real_+}$ to another function
 $c(h) \in \real^{\real_+}$.
 We also consider a mapping $\mathcal{M}$, called the \textit{mechanism}, 
 defined on $\mathcal{C}$ by matching a \textit{code} $c\in \mathcal{C}$
 to a set $\mathcal{M}(c)$ of code tuples.
\begin{definition} 
 The mechanism $\mathcal{M}$ is defined by 
 $  \mathcal{M}(  {\rm Id}_i ) = \{ f_i \}$ and 
\begin{equation}
  \label{jklfd243} 
\mathcal{M} ( g^* ) = \big\{ (f_1^*, (\partial_1 g)^*), (f_2^*,(\partial_2 g)^*), \ldots, (f_d^*,(\partial_d g)^*) \big\}, 
\end{equation}
 for $g$ a smooth function from $\real^d$ into $\real$.
\end{definition} 
In the next key lemma we show that $c(y)$ satisfies a system
of equations indexed by $c\in \mathcal{C}$.
 \begin{lemma}
  \label{l1}
  For any code $c\in \mathcal{C}$ we have
\begin{equation} 
 \label{s1}
 c(y) (t) = c(y)(0) + \sum_{Z \in \mathcal{M}(c)} \int_0^t \prod_{z \in Z} z(y)(s) ds,
 \qquad t\in \real_+. 
 \end{equation} 
\end{lemma} 
 \begin{Proof}
  When $c={\rm Id}_i$ we have
$$
 c(y)(t) = y_i(t) = y_i(0) + \int_0^t f_i(y_1(s),\ldots , y_d(s)) ds
 = y_i(0) + \int_0^t f_i^*(y_1,\ldots , y_d)(s) ds,  
$$
 hence \eqref{s1} holds since $\mathcal{M} ( {\rm Id}_i ) = \{ f^*_i \}$,
 $i=1,\ldots ,d$. 
 When $c = g^*\in \mathcal{C}$ with $c\not= {\rm Id}_i$, the equation 
$$
 g(y_1(t),\ldots , y_d(t)) = g(y_1(0),\ldots , y_d(0))
 + \sum_{j=1}^d \int_0^t f_j(y_1(s),\ldots , y_d(s)) \partial_j g(y_1(s),\ldots , y_d(s)) ds 
$$
 satisfied by $g^*(y)(t)$ reads 
\begin{equation}
\nonumber 
 g^*(y)(t) = g^*(y)(0) + \sum_{j=1}^d \int_0^t f_j^*(y)(s) ( \partial_j g)^*(y)(s) ds, 
\end{equation}
 and \eqref{s1} follows 
 by the definition \eqref{jklfd243} of $\mathcal{M}$. 
\end{Proof}
 We note that for any $g^* \in \mathcal{C}$ it is always possible to compute
 $g^*(y)(0)$ by applying the code $g^*$ to the solution $y$ of the ODE and then
 evaluating it at time $0$ as
 $g^*(y)(0) = g(y_0)$. In particular, the full knowledge of the function $y$
 is not necessary to compute $g^*(y)(0)$.

\subsubsection*{Example - One-dimensional autonomous ODE}
 Consider the solution $y(t)$ of the one-dimensional ODE 
\begin{equation}
\nonumber 
y(t) = y_0 + \int_0^t f(y(s)) ds, \qquad t\in \real_+, 
\end{equation}
 where $f(y(s))$ is expanded as 
 \begin{equation}
\nonumber 
 f(y(s)) = f(y_0) + \int_0^s f(y(u)) f'(y(u)) du, \qquad s\in \real_+.  
\end{equation}
 Here, $g^*$ denotes the mapping 
 \begin{align}
\nonumber 
    g^* : \real^{\real_+} & \longrightarrow \real^{\real_+}
    \\
\nonumber 
    (t\mapsto y(t)) & \longmapsto g^*(y):=(t \mapsto g(y(t))), 
\end{align} 
 for $g$ a smooth function from $\real$ to $\real$, 
 the set of codes is given by 
$$
 \mathcal{C}:= \left\{
         {\rm Id},
         \ \big( f^{(k)}\big)^*,
         ~~ k \geq 0 
         \right\}, 
$$
 where $\big( f^{(k)} \big)^*$, $k \geq 0$, denotes the operator acting as 
$$
\big( f^{(k)} \big)^* (y)(s) := f^{(k)} ( y(s)),
\qquad
 s \in \real_+, 
$$
 and the mechanism $\mathcal{M}$ is given by 
 $\mathcal{M} ( {\rm Id} ) = \{ f^* \}$ and 
 $\mathcal{M} ( g^* ) = \big\{ (f^*,(g')^*) \big\}$. 
\subsubsection*{Example - Non-autonomous ODE} 
 Consider the non-autonomous ODE 
\begin{equation}
\label{Eint3-1}
y(t) = y_0 + \int_0^t f_2 (s,y(s)) ds, \qquad t\in \real_+, 
\end{equation}
where $f_2$ is a smooth Lipschitz function defined on a domain of $\real_+ \times \real$.
This ODE can be rewritten as the system
\begin{equation}
\label{Eautsys1}
\left\{ 
\begin{array}{l}
 y_1(t) = t = \displaystyle y_1(0) + \int_0^t f_1(y_1(s), y_2(s)) ds\\
 y_2(t) = y(t) = \displaystyle y_2(0) + \int_0^t f_2(y_1(s), y_2(s)) ds
 \\
\end{array}
\right. 
\end{equation}
by taking $f_1\equiv 1$ and $y_1(s)=s$, $s\in \real_+$.
Here, the set of codes satisfies 
\begin{eqnarray*} 
\mathcal{C}:
 & = & \left\{
         {\rm Id}_1, {\rm Id}_2,
         \ \big(\partial_1^{i_1} \partial_2^{i_2} f_1\big)^*,
         \ \big(\partial_1^{j_1} \partial_2^{j_2} f_1\big)^*,
         ~~ i_1,i_2,j_1,j_2 \geq 0
         \right\}
         \\
          & = & \left\{
         {\rm Id}_1, {\rm Id}_2,
        0 , 1 , 
         \big(\partial_1^{j_1} \partial_2^{j_2} f_1\big)^*,
         ~~ j_1,j_2 \geq 0
         \right\}, 
\end{eqnarray*} 
where $0$ and $1$ denote constant functions. 
The mechanism $\mathcal{M}$ is given by 
$\mathcal{M} ( {\rm Id}_1 ) = \{ f^*_1 \} = \{ 1 \}$, 
$\mathcal{M} ( {\rm Id}_2 ) = \{ f^*_2 \}$, 
and 
$\mathcal{M} ( g^* ) = \big\{ ( f_1^* , (\partial_1 g)^*) ,
( f_2^* , (\partial_2 g)^*)
\big\}
= \big\{ ( 1 , (\partial_1 g)^*) ,
( f_2^* , (\partial_2 g)^*)
 \big\}$ for $g$ a smooth function from $\real^2$ into $\real$. 

 \medskip

 More generally, any non autonomous system can be transformed
into an autonomous system by addition of a dimension.
In particular, any higher order ordinary differential equation of the form
$$
y^{(d)}(t) = f\big(t,y(t),y'(t),\ldots , y^{(d-1)}(t)\big)
$$
can be written as a system of the form \eqref{Eautsys} by
taking $f_1(y_1,\ldots , y_d) \equiv 1$, 
$f_i(y_1,\ldots , y_d) := y_{i+1}$, $i=2,\ldots , d-1$ and
$f_d(y_1,\ldots , y_d) := f(y_1,\ldots ,y_d)$,
with
$y_1(t)=t$ and $y_i(t)=y^{(i-1)}(t)$, $i=2,\ldots , d$.

\section{Coding trees} 
 For each code $c \in \mathcal{C}$ we denote by $I_c$ a uniformly distributed 
random variable on $\mathcal{M}(c)$. 
For example, 
when $c = g^*$,
 since $\mathcal{M} ( g^* ) = \big\{ (f_1^*, (\partial_1 g)^*), (f_2^*,(\partial_2 g)^*), \ldots, (f_d^*,(\partial_d g)^*) \big\}$,
we have 
$$
\mathbb{P}(I_{g^*} = (f_i^*, (\partial_i g)^*)
\big) = \frac{1}{d}, \qquad i = 1,\ldots, d.
$$
 In the sequel we will use the notation 
 $$
 q_c(b) :=
 \mathbb{P}(I_c=b), \qquad
 b\in \mathcal{M}(c), \quad c \in \mathcal{C}.
$$
 In addition, we consider
 a probability density function $\rho: \real_+ \to (0,\infty )$
 and 
 \begin{itemize}
\item an i.i.d. family $(\tau^{i,j})_{i,j\geq 1}$ of random variables
 with distribution $\rho (t)dt$ on $\real_+$,
\item for each $c\in \mathcal{C}$, an i.i.d. family $(I_c^{i,j})_{i,j\geq 1}$ of discrete
  random variables, with
  $$
  \mathbb{P}\big( I_c^{i,j}=b \big) = q_c(b) >0,
  \qquad b \in \mathcal{M}(c),
  $$
\end{itemize}
 where the sequences $(\tau^{i,j})_{i,j\geq 1}$ and $(I_c^{i,j})_{c\in \mathcal{C}, i,j\geq 1}$ are assumed to be mutually independent.
 For every $n \geq 1$ we also consider an
 injection $\pi_n : \mathbb{N}^n \rightarrow \mathbb{N}$.

 \medskip

 Let $t>0$.
 For each $i\in \{1,\ldots , d\}$ we construct a random tree
 starting from an initial particle labelled $\widebar{1} := (1)$
 bearing the code ${\rm Id}_i$
 at time $0$, which lives up to a random time $\tau^{1,1}$
distributed according to $\rho$. If $\tau^{1,1} > t$, the branching process stops.
 Otherwise, if $\tau^{1,1} \leq t$, a new particle with label $(1,1)$ is created,
 and bears the code $f_i^*$ since $\mathcal{M} ( {\rm Id}_i ) = \{ f_i^* \}$, 
 and independently follows the same pattern as the first one. 
 This new branch lives during the time $\tau^{2,\pi_2{(1,1)}}$. 
 If $ \tau^{(1,1)} + \tau^{2,\pi_2{(1,1)}} > t$ then the tree stops branching,
 otherwise, if $\tau^{1,1} + \tau^{2,\pi_2(1,1)} \leq t $,
 the particle branches in two branches
 $(f_i^*,(\partial_j f_i)^* )$ 
 chosen uniformly in 
 $\mathcal{M} ( f_i^* ) = \big\{ (f_1^*,(\partial_1 f_i)^* ), (f_2^*,(\partial_2 f_i)^* ), \ldots, (f_d^*,(\partial_d f_i)^* ) \big\}$.

 \medskip
 
More generally, a particle with code $c\in \mathcal{C}$
at the generation $n\geq 1$ is assigned a label of the form
$\bar{k} = (1,k_2,\ldots ,k_n) \in \mathbb{N}^n$,
 while its parent label is 
$\bar{k}{\text{--}} := (1,k_2,\ldots ,k_{n-1})$.
The birth time
of particle $\bar{k}$ is denoted by $T_{\bar{k}{\text{--}}}$, 
and its lifetime $\tau^{n,\pi_n(\bar{k})}$ is the element of index
$j=\pi_n(\bar{k})$ in the i.i.d. sequence
$(\tau^{n,j})_{j\geq 1}$.
 If $T_{\bar{k}} := \T_{\bar{k}{\text{--}}} + \tau^{n,\pi_n(\bar{k})} < t$,
 we draw a sample
 $I_c^{n,\pi_n(\bar{k})} = (c_1,\ldots ,c_l)$ uniformly in $\mathcal{M}(c)$
 with $l\in \{1,2\}$,
 and the particle $\bar{k}$ branches into
 $\big|I_c^{n,\pi_n(\bar{k})}\big|\in \{1,2\}$ offsprings at generation $(n+1)$,
 which are labeled by $\bar{k}=(1,\ldots ,k_n,j)$, $j=1,\ldots ,\big|I_c^{n,\pi_n(\bar{k})}\big|$.
 The particle with label ending with an integer $j$ will carry the code $c_j$. 
 Finally, the code of particle $\bar{k}$ will be denoted by
 $c_{\bar{k}} \in \mathcal{C}$.
 The death time of the particle $\widebar{k}$ is 
 $T_{\widebar{k}}$, and its birth time is $T_{\widebar{k}{\text{--}}}$.
 
\begin{definition}
  We denote by $\mathcal{T}_{t,c}$ the random tree
  constructed from the above argument
  started from any code $c\in \mathcal{C}$. 
\end{definition}
The family $(\mathcal{T}_{t,{\rm Id}_i })_{i=1,\ldots , d}$
 of trees will be used for the stochastic representation of the solution
 $y(t)=(y_1(t),\ldots , y_d(t))$ of the ODE system \eqref{Eautsys}, while the trees $\mathcal{T}_{t,c}$ 
 will be used for the stochastic representation of $c(y)(t)$. 
 The next table summarizes the notation introduced so far.

 \medskip

\begin{center}
  \centering
  \begin{tabular}{||l | c||}
 \hline
 Object & Notation \\ [0.5ex]
 \hline\hline
 Initial time  & $0$ \\
 \hline
 Tree ending at time $t$ with initial code $c$ &  $\mathcal{T}_{t,c}$ \\  \hline
 Particle (or label) of generation $n\geq 1$ & $\widebar{k}=(1,k_2,\ldots ,k_n)$\\
 \hline
 First branching time & $T_{\widebar{1}}$\\
 \hline
Birth time of a particle $\widebar{k}$ & $T_{\widebar{k}{\text{--}}}$ \\
\hline
Death time of a particle $\widebar{k}$ & $T_{\widebar{k}}$ \\
 \hline
 Lifespan of a particle & $ 
 T_{\widebar{k}} - T_{\widebar{k}\text{--}}$ \\
\hline
Code of a particle $\widebar{k}$ & $c_{\widebar{k}}$
\\ 
\hline
\end{tabular}
\end{center}

\vspace{0.05cm}

\noindent
 The following graph represents a sample of the random tree $\mathcal{T}_{t,{\rm Id}_i}$, $i=1,\ldots , d$. 

\tikzstyle{level 1}=[level distance=4cm, sibling distance=4cm]
\tikzstyle{level 2}=[level distance=5cm, sibling distance=3cm]

\begin{center}
\resizebox{0.94\textwidth}{!}{
\begin{tikzpicture}[scale=1.0,grow=right, sloped]
\node[rectangle,draw,black,text=black,thick]{$0$}
    child {
        node[rectangle,draw,black,text=black,thick] {$T_{\widebar{1}}$}
            child {
                node[rectangle,draw,black,text=black,thick] {$T_{(1,1)}$}
                child{
                node[rectangle,draw,black,text=black,thick,yshift=-0.4cm]{$T_{(1,1,2)}$}
                    child{
                    node[rectangle,draw,black,text=black,thick]{$t$}
                    edge from parent
                    node[above]{$(1,1,2,2)$}
                    node[below]{$(\partial_k\partial_j f_i)^*$}
                    }
                    child{
                    node[rectangle,draw,black,text=black,thick]{$t$}
                    edge from parent
                    node[above]{$(1,1,2,1)$}
                    node[below]{$f^*_k$}
                    }
                edge from parent
                node[above]{$(1,1,2)$}
                node[below]{$(\partial_j f_i)^*$}
                }
                child{
                node[rectangle,draw,black,text=black,thick,yshift=0.4cm]{$T_{(1,1,1)}$}
                    child{
                    node[rectangle,draw,black,text=black,thick]{$t$}
                    edge from parent
                    node[above]{$(1,1,1,2)$}
                    node[below]{$(\partial_l f_j)^*$}
                    }
                    child{
                    node[rectangle,draw,black,text=black,thick]{$t$}
                    edge from parent
                    node[above]{$(1,1,1,1)$}
                    node[below]{$f^*_l$}
                    }
                edge from parent
                node[above]{$(1,1,1)$}
                node[below]{$f^*_j$}
                }
                edge from parent
                node[above] {$(1,1)$}
                node[below]  {$f^*_i$}
            }
            edge from parent
            node[above] {$\widebar{1}$}
            node[below]{${\rm Id}_i$}
    };
\end{tikzpicture}
}
\end{center}

\section{Probabilistic representation of ODE solutions}
\label{s3}
 In this section we define the random multiplicative functional which will be used
 to represent ODE solutions.
 We let $\widebar{F}$ denote the tail distribution function of
 $\rho$, i.e. 
 $$
 \widebar{F}(t) := \int_t^\infty \rho(u) du, \qquad t \in \real_+. 
 $$
 We denote by $\mathcal{K}^\circ$ the set of particles that do not live until time $t$,
 and by $\mathcal{K}^{\partial}$ the set of particles that die after time $t$.
 \begin{definition}
   Given $\mathcal{T}_{t,c}$ a random coding tree
   started with the code $c \in \mathcal{C}$ 
  and ending at time $t\geq 0$, 
    we define the universal multiplicative functional $\mathcal{H}$ by 
$$
{\mathcal{H}(\mathcal{T}_{t,c}) := \prod_{\widebar{k} \in \mathcal{K}^{\circ}} \frac{1}{q_{c_{\widebar{k}}}(I_{c_{\widebar{k}}})\rho(T_{\widebar{k}}-T_{\widebar{k}{\text{--}}})} \prod_{\widebar{k} \in \mathcal{K}^{\partial}} \frac{c_{\widebar{k}}(y)\big(0\big)}{\widebar{F}(t-T_{\widebar{k}{\text{--}}})}}.
$$
\end{definition}{}
 The next result gives the probabilistic representation of ODE solutions
 as an expected value over random coding trees.
 \begin{theorem}
   \label{t1}
   Let $T>0$ for which there exists $K_0>0$ such that
   $$
   \E \big[ \big| \mathcal{H}(\mathcal{T}_{t,c}) \big| \big]
\leq K_0,
\qquad c\in \mathcal{C},
\quad
t\in [0,T]. 
$$
Then, for any $c\in \mathcal{C}$ we have the probabilistic representation
 \begin{equation}
   \label{fjhkldsf} 
   c(y)(t) = \E \big[ \mathcal{H}(\mathcal{T}_{t,c}) \big], \qquad
   t\in [0,T],
 \end{equation}
 where $y(t)= (y_1(t),\ldots , y_d(t))$ is the solution of the system of ODEs
 \eqref{Eautsys}.
 In particular, taking $c={\rm Id}_i$, we have 
 \begin{equation}
 \nonumber 
 y_i (t)
 =
 y_{{\rm Id}_i} (t) = \E \big[ \mathcal{H}(\mathcal{T}_{t,{\rm Id}_i}) \big],
    \qquad
    t\in [0,T],
    \quad i=1,\ldots , d.
\end{equation}
\end{theorem}
\begin{Proof}
   For $c\in \mathcal{C}$ we let 
$$
   y_c(t) := \E \big[ \mathcal{H}(\mathcal{T}_{t,c})\big],
   \qquad t\in [0,T]. 
$$
By conditioning on the first branching time $T_{\widebar{1}}$, 
for all $i=1,\ldots , d$ the first particle bearing the code ${\rm Id}_i$
branches at time $T_{\widebar{1}}$
into a new particle bearing the code $f^*_i$
as $\mathcal{M} ( {\rm Id}_i ) = \{ f^*_i \}$, hence we have 
\begin{equation*}
\begin{split}
y_{{\rm Id}_i}(t) &= \E \big[ \mathcal{H}(\mathcal{T}_{t,{\rm Id}_i })\mathbbm{1}_{\{ T_{\widebar{1}}> t \}} + \mathcal{H}(\mathcal{T}_{t,{\rm Id}_i }) \mathbbm{1}_{\{ T_{\widebar{1}}\leq t \}} \big]\\
&= \E \Bigg[ \frac{y_i(0)}{\widebar{F}(t)} \mathbbm{1}_{\{ T_{\widebar{1}}>t \}} \Bigg]
+  \E \Bigg[
  \frac{y_{f^*_i}(t-T_{\widebar{1}})}{\rho(T_{\widebar{1}})}\mathbbm{1}_{\{ T_{\widebar{1}} \leq t \}}  \Bigg]\\
&= y_i(0) \frac{\mathbb{P}(T_{\widebar{1}}>t)}{ \widebar{F}(t)} + \int_0^t \frac{y_{f^*_i}(t-s)}{\rho(s)} \rho(s) ds\\
&= y_i(0) + \int_0^t y_{f^*_i}(s) ds,
   \qquad t\in [0,T]. 
\end{split}
\end{equation*}
 Similarly, starting from any code $g^* \in \mathcal{C}$ different from ${\rm Id}_i$,
 the particle branches at time $T_{\widebar{1}}$ in $d$ possible different ways 
 into two particles with codes $(f^*_j,(\partial_j g)^*)$, $j=1,\ldots , d$, 
 hence we have 
\begin{align}
\nonumber
 y_{g^*} (t) & = \E \big[ \mathcal{H}(\mathcal{T}_{t,{g^*}})\mathbbm{1}_{\{ T_{\widebar{1}}> t \}} + \mathcal{H}(\mathcal{T}_{t,{g^*}}) \mathbbm{1}_{\{ T_{\widebar{1}}\leq t \}} \big]
\\
\nonumber
&= \E \bigg[ \frac{{g^*}(y_0)}{\widebar{F}(t)} \mathbbm{1}_{\{ T_{\widebar{1}}>t \}}
  +
  \mathbbm{1}_{\{ T_{\widebar{1}} \leq t \}}
  \sum_{j=1}^d
  {\bf 1}_{\{ I_{g^*} = (f^*_j,(\partial_j g)^*)\}}
 \frac{y_{f^*_j}(t-T_{\widebar{1}})y_{(\partial_j g)^*}(t-T_{\widebar{1}})}{  q_{g^*}(I_{g^*}) \rho(T_{\widebar{1}})}
\bigg]
\\
\nonumber
&= \E \bigg[ \frac{{g^*}(y_0)}{\widebar{F}(t)} \mathbbm{1}_{\{ T_{\widebar{1}}>t \}}
  +
  \mathbbm{1}_{\{ T_{\widebar{1}} \leq t \}}
\sum_{j=1}^d 
q_{g^*_j}((f^*_j,(\partial_j g)^*)) 
\frac{y_{f^*_j}(t-T_{\widebar{1}})y_{(\partial_j g)^*}(t-T_{\widebar{1}})}{ q_{g^*_j}((f^*_j,(\partial_j g)^*)) \rho(T_{\widebar{1}})}
\bigg]
\\
\nonumber
&= \E \bigg[ \frac{{g^*}(y_0)}{\widebar{F}(t)} \mathbbm{1}_{\{ T_{\widebar{1}}>t \}}
  \bigg]
+
\sum_{Z\in {\cal M}({g^*})}
\E \bigg[  \mathbbm{1}_{\{ T_{\widebar{1}} \leq t \}}
 \frac{1}{  \rho(T_{\widebar{1}})}
 \prod_{z\in Z}
  y_z(t-T_{\widebar{1}})
\bigg]
\\
\nonumber
& = {g^*}(y)(0) + \sum_{Z \in \mathcal{M}({g^*})} \int_0^t \frac{1}{\rho(s)}
\left( \prod_{z \in Z} y_z(t-s) \right) \rho (s) ds, 
\\
  \label{fjhkdsf} 
& = {g^*}(y)(0) + \sum_{Z \in \mathcal{M}({g^*})} \int_0^t \prod_{z \in Z} y_z(s) ds, 
\end{align}
 which yields the system of equations
\begin{equation} 
  \label{s}
  y_c (t) = c(y)(0) + \sum_{Z \in \mathcal{M}(c)} \int_0^t \prod_{z \in Z} y_z(s) ds,    \qquad t\in [0,T], \quad c \in \mathcal{C}. 
 \end{equation} 
By the Cauchy-Lipschitz theorem on the Banach space of sequences 
$\ell^\infty$, 
this system admits a unique maximal solution. 
 We conclude by noting that from Lemma~\ref{l1}, the family of functions  
 $(c (y))_{c\in \mathcal{C} }$ is the solution of the system \eqref{s}, 
 hence
 $(c (y))_{c\in \mathcal{C} } = (y_c)_{c\in \mathcal{C} }$, and 
$$
\E \big[ \mathcal{H}(\mathcal{T}_{t,c}) \big]
=
y_c (t) = 
c(y)(t), \qquad t\in [0,T]. 
$$
\end{Proof} 
In numerical applications
the expected value $\E \big[ \mathcal{H}(\mathcal{T}_{t,c}) \big]$
in Theorem~\ref{t1} is estimated as the average
 $$
 \frac{1}{N}\sum_{k=1}^N \mathcal{H}(\mathcal{T}_{t,c})^{(k)}
 $$
 where
 $\mathcal{H}(\mathcal{T}_{t,c})^{(1)}, \ldots , \mathcal{H}(\mathcal{T}_{t,c})^{(N)}$
 are independent samples of $\mathcal{H}(\mathcal{T}_{t,c})$. 
 In this case, the error on the estimate of
 $\E \big[ \mathcal{H}(\mathcal{T}_{t,c}) \big]$
 from the Monte Carlo method can be estimated as the standard deviation
$$
\left( \E \left[ \left( \E \big[ \mathcal{H}(\mathcal{T}_{t,c}) \big] 
   - \frac{1}{N}\sum_{k=1}^N \mathcal{H}(\mathcal{T}_{t,c})^{(k)}
   \right)^2 \right]
\right)^{1/2} =
 \frac{1}{\sqrt{N}} \sqrt{\Var \big[ \mathcal{H}(\mathcal{T}_{t,c}) \big]}. 
$$
\subsubsection*{Complexity of the algorithm} 
We note that the complexity of the algorithm grows linearly
with the dimension $d$, as $d$ trees
$({\cal T}_{t,{\rm Id}_i})_{i=1,\ldots , d}$ are used to
generate the multidimensional solution $(y_1,\ldots , y_d)$,
whereas the complexity of finite difference methods is generally polynomial
in the dimension $d$, depending on the order chosen in the truncation of
\eqref{sdfh}. 
On the other hand, all
trees are at most binary regardless of the dimension $d$, as shown by the construction
\eqref{jklfd243} of the mechanism $\mathcal{M}$.

\medskip

Regarding complexity in time, let $\ell_c(t)$ denote the mean size
of the random tree ${\cal T}_{t,c}$ generated until time $t\geq 0$,
with, by construction, 
$\ell (t) := \ell_{{\rm Id}_i} (t)$ for all $i=1,\ldots , d$,
and 
$m (t) := \ell_c(t)$ for all $c\notin \{{\rm Id}_i\}_{i=1,\ldots , d}$. 
The same argument as in the proof of
Theorem~\ref{t1} shows that $(\ell (t),m(t))$ satisfies the sytem of
 integral equations 
\begin{equation}
\nonumber 
\left\{
\begin{array}{l}
  \displaystyle
  \ell (t) = \int_t^\infty \rho (s) ds + \int_0^t \rho (s) m (t-s) ds, 
  \medskip
  \\
 \displaystyle
  m (t) = \int_t^\infty \rho (s) ds + 2 \int_0^t \rho (s) m (t-s) ds. 
\end{array}
\right.
\end{equation} 
 When $\tau$ has the exponential density $\rho (s) = \lambda \re^{- \lambda s}$
 with parameter $\lambda > 0$, this leads to the system 
\begin{equation}
\nonumber 
\left\{
\begin{array}{l}
  \displaystyle
  \ell' (t) = - \lambda \re^{-\lambda t} - \lambda^2 \int_0^t  \re^{-\lambda (t-s)} m (s) ds + \lambda m (t) = \lambda ( m (t) - \ell (t ) ) 
  \medskip
  \\
 \displaystyle
 m' (t) = - \lambda \re^{-\lambda t} - 2 \lambda^2 \int_0^t  \re^{-\lambda (t-s)} m (s) ds + 2  \lambda m (t)
 = \lambda m (t),  
\end{array}
\right.
\end{equation} 
 with solution 
$$
 \ell (t) = \cosh ( \lambda t), \quad
 m(t) = \re^{\lambda t} = \E \big[ 2^{N_t} \big], \qquad t\geq 0,
$$
 and
 where $(N_t)_{t\in \real_+}$ is a standard Poisson process with intensity
 $\lambda > 0$.
 We note that this estimate remains the same independently of dimension $d\geq 1$. 
\subsubsection*{Integrability condition}
 The following proposition provides sufficient conditions
 ensuring that the representation formula \eqref{fjhkldsf} of
 Theorem~\ref{t1} holds at any time within certain interval.  
 In order to represent an ODE solution beyond that time interval we may
 reuse the numerical value obtained close to its boundary
 as new initial condition in order to represent the solution on an extended
 time interval. 
 We note that the constant $K$ in the next proposition
 depends on $y_0$ and $f$ and its derivatives,
 and is independent of the system dimension $d \geq 1$. 
 However, the integrability condition \eqref{fjkld} 
 is stronger in higher dimensions.   
\begin{prop}
\label{p1}
 Assume that there exists $K>0$ such that $c(y)(0) \leq K$
 for any $c \in \mathcal{C}$,
 that the density function $\rho$ is nonincreasing, and that 
  \begin{equation}
    \label{fjkld} 
\rho(T) \geq d, \quad K \leq \widebar{F}(T). 
\end{equation} 
 Then, there exists $K(T)>0$ such that 
\begin{equation}
  \label{gfdg} 
   \E \big[ \big| \mathcal{H}(\mathcal{T}_{s,c})\big| \big] 
\leq K(T),
\quad c\in \mathcal{C},
\quad
t\in [0,T]. 
\end{equation} 
\end{prop} 
 \begin{Proof}
 Under Condition~\eqref{fjkld}, 
 since $q_{\min} := \min_{c\in \mathcal{C}} q_c(I_c) = 1/d$, we have 
$$
 \prod_{\widebar{k} \in \mathcal{K}^{\circ}} \frac{1}{q_{\min}\rho(T_{\widebar{k}}-T_{\widebar{k}{\text{--}}})} \prod_{\widebar{k} \in \mathcal{K}^{\partial}} \frac{K }{\widebar{F}(t-T_{\widebar{k}{\text{--}}})} \leq 1, 
$$
 hence 
$$
 \E \big[ \big| \mathcal{H}(\mathcal{T}_{t,c})\big| \big] \leq
 \E_c \left[ \prod_{\widebar{k} \in \mathcal{K}^{\circ}} \frac{1}{q_{\min}\rho(T_{\widebar{k}}-T_{\widebar{k}{\text{--}}})} \prod_{\widebar{k} \in \mathcal{K}^{\partial}} \frac{K }{\widebar{F}(t-T_{\widebar{k}{\text{--}}})} \right] \leq 1,
 \qquad
 t\in [0,T]. 
$$
\end{Proof}
 We note that trying to relax the integrability condition
 \eqref{fjkld} by choosing a higher $\rho (T)$
 will result into a smaller value of $\E [ \tau ]$,
 therefore increasing the time complexity of the algorithm. 
 On the other hand, assuming that 
$$
  S_c(t): =
 \E_c \Bigg[ \prod_{\widebar{k} \in \mathcal{K}^{\circ}} \frac{1}{q_{c_{\widebar{k}}}(I_{c_{\widebar{k}}})\rho(
     T_{\widebar{k}}-T_{\widebar{k}{\text{--}}}
          )} \prod_{\widebar{k} \in \mathcal{K}^{\partial}} \frac{K }{\widebar{F}(t-T_{\widebar{k}{\text{--}}})} \Bigg] < \infty, 
   $$
  yields the system of equations 
\begin{equation*}
\begin{cases}
  \displaystyle
  S_{{\rm Id}_i}(t) = K + \int_0^t S_{f^*_i}(s) ds
  \\
  \\
\displaystyle
S_{g^*} (t) = K + \sum_{i=1}^d \int_0^t S_{f^*_i}(s) S_{(\partial_i g)^*}(s) ds. 
\end{cases}
\end{equation*}
   Since the corresponding trees
   ${\cal T}_{t,f_i^*}$, ${\cal T}_{t,c}$ 
   have same random shape distribution, we have 
   $S_c(t) = S_{f^*_i}(t)$ for all codes $c=g^* \in {\cal C}$ such that 
   $g = \partial_1^{i_1}\cdots \partial_d^{i_d} f_i$,
   $i_1,\ldots , i_d \geq 0$, $i =1, \ldots, d$.
   This gives the system 
\begin{equation*}
\begin{cases}
  \displaystyle
  S_{{\rm Id}_i}(t) = K + \int_0^t S_{f^*_i}(s) ds
  \\
  \\
\displaystyle
S_{f^*_i} (t) = K + d \times \int_0^t ( S_{f^*_i}(s) )^2 ds,
\end{cases}
\end{equation*}
which can be solved as 
$$
 S_{{\rm Id}_i}(t) = K - \frac{1}{d} \log(1-Ktd), 
 \qquad 
 S_{f^*_i}(t) = \frac{K}{1-Ktd}, 
 \qquad
  i = 1,\ldots , d, 
$$
 hence the finiteness of $S_c(t)$ holds at most until time 
\begin{equation}
   \label{p1.1} 
  T < \frac{1}{Kd}. 
\end{equation}
  
\begin{remark} 
 Although the non-autonomous ODE \eqref{Eint3-1}
 can be treated using an
 autonomous $2$-dimensional system of the form \eqref{Eautsys1},
 the probabilistic representation \eqref{fjhkldsf} of its solution 
 can also be obtained using a single random tree.
 For this, we expand $f_2(s,y(s))$ as 
 \begin{equation}
\nonumber 
 f_2(s,y(s)) = f_2(0,y_0) + \int_0^s \big(
 \partial_0 f_2(u,y(u))
 + f_2(u,y(u)) \partial_1 f_2(u,y(u))
 \big) du. 
\end{equation}
 In this case, the set of codes is defined as 
$$
 \mathcal{C}:= \left\{
         {\rm Id},
         \ \big(\partial_0^k \partial_1^l f_2 \big)^*,
         ~~ k,l \geq 0 
         \right\}, 
$$
         and the mechanism $\mathcal{M}$ is given by
         $$\mathcal{M} ( {\rm Id} ) = \{ f_2^* \}
         \quad
         \mbox{and}  
         \quad
         \mathcal{M} ( g^* ) = \big\{ (\partial_0 g)^* , (f_2^*,(\partial_1 g)^*) \big\}.
         $$
 A sample of this random tree is presented below. 
 
\tikzstyle{level 1}=[level distance=4cm, sibling distance=4cm]
\tikzstyle{level 2}=[level distance=5cm, sibling distance=3cm]

\begin{center}
\resizebox{0.94\textwidth}{!}{
\begin{tikzpicture}[scale=1.0,grow=right, sloped]
\node[rectangle,draw,black,text=black,thick]{$0$}
    child {
        node[rectangle,draw,black,text=black,thick] {$T_{\widebar{1}}$}
            child {
                node[rectangle,draw,black,text=black,thick] {$T_{(1,1)}$}
                child{
                node[rectangle,draw,black,text=black,thick]{$T_{(1,1,2)}$}
                    child{
                    node[rectangle,draw,black,text=black,thick]{$t$}
                    edge from parent
                    node[above]{$(1,1,2,2)$}
                    node[below]{$(\partial_1^2 f_2)^*$}
                    }
                    child{
                    node[rectangle,draw,black,text=black,thick]{$t$}
                    edge from parent
                    node[above]{$(1,1,2,1)$}
                    node[below]{$f_2^*$}
                    }
                edge from parent
                node[above]{$(1,1,2)$}
                node[below]{$(\partial_1 f_2)^*$}
                }
                child{
                node[rectangle,draw,thick]{$T_{(1,1,1)}$}
                    child{
                    node[rectangle,draw,black,text=black,thick]{$t$}
                    edge from parent
                    node[above]{$(1,1,1,1)$}
                    node[below]{$(\partial_0 f_2)^*$}
                    }
                edge from parent
                node[above]{$(1,1,1)$}
                node[below]{$f_2^*$}
                }
                edge from parent
                node[above] {$(1,1)$}
                node[below]  {$f_2^*$}
            }
            edge from parent
            node[above] {$\widebar{1}$}
            node[below]{${\rm Id}$}
    };
\end{tikzpicture}
}
\end{center}
 In this case, the system of equations satisfied by $S_c(t)$ can be written as 
\begin{equation*}
\begin{cases}
  \displaystyle
  S_{{\rm Id}}(t) = K + \int_0^t S_{f^*_i}(s) ds
  \\
  \\
\displaystyle
  S_{f^*_i} (t) = K + \int_0^t (S_{f^*_i}(s))^2 ds + \int_0^t S_{f^*_i}(s)ds, 
\end{cases}
\end{equation*}
 which can be solved as 
$$
 S_{{\rm Id}}(t) = K - \log(1+K(1-e^t)), 
 \qquad 
 S_{f^*_i}(t) = \frac{Ke^t}{1+K(1-e^t)},
 \qquad
  i = 1,\ldots , d.
$$
 Therefore, the finiteness of $S_c(t)$ 
 holds until time 
 \begin{equation}
   \label{fn3} 
 T < \log\left( 1 + \frac{1}{K}\right). 
\end{equation}
Under the conditions $\rho(T) \geq 2$ and $K \leq \widebar{F}(T)$, 
since $q_{\min} := \min_{c\in \mathcal{C}} q_c(I_c) = 1/2$, we have 
$$
 \prod_{\widebar{k} \in \mathcal{K}^{\circ}} \frac{1}{q_{\min}\rho(T_{\widebar{k}}-T_{\widebar{k}{\text{--}}})} \prod_{\widebar{k} \in \mathcal{K}^{\partial}} \frac{K }{\widebar{F}(t-T_{\widebar{k}{\text{--}}})} \leq 1, 
$$
 hence 
$$
 \E \big[ \big| \mathcal{H}(\mathcal{T}_{t,c})\big| \big] \leq
 \E_c \left[ \prod_{\widebar{k} \in \mathcal{K}^{\circ}} \frac{1}{q_{\min}\rho(T_{\widebar{k}}-T_{\widebar{k}{\text{--}}})} \prod_{\widebar{k} \in \mathcal{K}^{\partial}} \frac{K }{\widebar{F}(t-T_{\widebar{k}{\text{--}}})} \right] \leq 1,
 \quad t\in [0,T]. 
$$
 In this setting the integrability condition \eqref{fn3}
 is stronger than \eqref{p1.1} when $K$ is sufficiently large,
 and weaker otherwise.
 On the other hand,
 in comparison with the autonomous system \eqref{Eautsys1},
 time complexity of the algorithm is 
 divided by at least two due the use of a single coding tree with less branches. 
\end{remark} 
\section{Examples}
\label{s4}
\subsubsection*{Exponential series} 
\noindent
 We first consider the equation   
\begin{equation}
\nonumber 
\begin{cases}
y'(t) = y(t) \\
y(0)= y_0
\end{cases}
\end{equation} 
 rewritten in integral form as
\begin{equation}
\nonumber 
y(t) = y_0 + \int_0^t y(s)ds, \quad t\in \real_+, 
\end{equation}
 whose solution admits the power series expansion 
$$
y(t) = y_0 e^t = y_0 \sum_{n=0}^\infty \frac{t^n}{n!}, \qquad t\in \real_+.
$$
 Here we have $\mathcal{C}= \{ {\rm Id} \}$, 
 and the mechanism $\mathcal{M}$ satisfies $\mathcal{M} ( {\rm Id} ) = \{ {\rm Id} \}$. 
The particle of generation $n$ bears the label $\widebar{k} = (1,\ldots,1) \in \mathbb{N}^n$, and its parent is the particle $\widebar{k}{\text{--}} = (1,\ldots,1) \in \mathbb{N}^{n-1}$.
When the random times $(\tau^k)_{k\geq 1}$ are independent and exponentially distributed, 
i.e. $\rho(s) = e^{-s}$ with $\widebar{F}(t) = e^{-t}$, the total number of branches in the random tree $\mathcal{T}_t$ is given by $N_{t}+1$ where $(N_t)_{t \geq 0}$ is a standard Poisson process with unit intensity.

  \begin{center}
    \resizebox{0.9\textwidth}{!}
            {
\begin{tikzpicture}[scale=1,grow=right, sloped]
\node[rectangle,draw,black,text=black,thick]{$0$}
    child {
        node[rectangle,draw,black,text=black,thick] {$T_{\widebar{1}}$}
            child {
                node[rectangle,draw,black,text=black,thick] {$T_{(1,1)}$}
                child{
                    node[rectangle,draw,black,text=black,thick] {$T_{(1,1,1)}$}
                    child{
                        node[rectangle,draw,black,text=black,thick] {$t$}
                    edge from parent
                    node[above] {$(1,1,1,1)$}
                    node[below]  {${\rm Id}$}    
                    }
                edge from parent
                node[above] {$(1,1,1)$}
                node[below]  {${\rm Id}$}             
                }
                edge from parent
                node[above] {$(1,1)$}
                node[below]  {${\rm Id}$}
            }
            edge from parent
            node[above] {$\widebar{1}$}
            node[below]{${\rm Id}$}
    };
\end{tikzpicture}
            }
              \end{center}
            In this case, the multiplicative functional 
$$
\mathcal{H}(\mathcal{T}_t):= \prod_{k \in \mathcal{K}^\circ} \frac{1}{\rho(T_{\widebar{k}} - T_{\widebar{k}{\text{--}}})} \prod_{k \in \mathcal{K}^\partial} \frac{y_0}{\widebar{F}(t-T_{\widebar{k}{\text{--}}})}
$$
 simplifies to the deterministic expression 
$$
 \mathcal{H}(\mathcal{T}_t) =
 \left( \prod_{k=1}^{N_{t}} \frac{1}{e^{-\tau^k}}
 \right) 
 \frac{y_0}{e^{-(t-\sum_{k=1}^{N_{t}} \tau^k)}} 
= y_0 e^{t}, \qquad t\in \real_+, 
$$ 
 in which we take $\sum_{n=1}^0 1 := 0$ and $\prod_{n=1}^0 1:= 1$.

\subsubsection*{One-dimensional autonomous ODE} 
  \noindent
 Consider the autonomous ODE 
\begin{equation}
\label{Eint2}
y(t) = y_0 + \int_0^t f(y(s)) ds, \qquad t\in \real_+, 
\end{equation}
where $f \in C^\infty(\real;\real)$ is bounded together
with its derivatives $f^{(k)}$ of order $k\geq 1$, with
\begin{equation}
\nonumber 
  \big|f^{(k)}(y_0)\big| \leq K, \qquad k\geq 0,
  \end{equation}
  for some $K>0$. 
 For any $k\geq 0$ we have the integral equation 
\begin{equation*}
f^{(k)}(y(t)) = f^{(k)}(y_0) + \int_0^t f(y(s)) f^{(k+1)}(y(s)) ds, \qquad t\in \real_+, 
\end{equation*}
 with the set of codes 
$
 \mathcal{C} := \big\{ {\rm Id}, \ \big( f^{(k)} \big)^*, \ k\geq 0 \big\}$,
 and the mechanism $\mathcal{M}$ is defined by 
$$
\mathcal{M}({\rm Id}) := \{f^*\}, \quad 
\mathcal{M}(f^*):= \big\{(f^*,(f')^*)\big\}, \quad 
\mathcal{M} \big( \big(f^{(k)}\big)^* \big) := \big\{\big(f^*,\big(f^{(k+1)}\big)^* \big)\big\},
\quad k \geq 1. 
$$
 Below is a representation of a sample of the random tree $\mathcal{T}_{t,{\rm Id}}$.

\tikzstyle{level 1}=[level distance=3cm, sibling distance=1cm]
\tikzstyle{level 2}=[level distance=4cm, sibling distance=4cm]

\begin{center}
\resizebox{0.75\textwidth}{!}{
\begin{tikzpicture}[scale=1.0,grow=right, sloped]
\node[rectangle,draw,black,text=black,thick]{$0$}
    child {
        node[rectangle,draw,black,text=black,thick] {$T_{\widebar{1}}$}
            child {
                node[rectangle,draw,black,text=black,thick] {$T_{(1,1)}$}
                child{
                node[rectangle,draw,black,text=black,thick]{$T_{(1,1,2)}$}
                    child{
                    node[rectangle,draw,black,text=black,thick]{$t$}
                    edge from parent
                    node[above]{$(1,1,2,2)$}
                    node[below]{$(f'')^*$}
                    }
                    child{
                    node[rectangle,draw,black,text=black,thick]{$t$}
                    edge from parent
                    node[above]{$(1,1,2,1)$}
                    node[below]{$f^*$}
                    }
                edge from parent
                node[above]{$(1,1,2)$}
                node[below]{$(f')^*$}
                }
                child{
                node[rectangle,draw,thick]{$t$}
                edge from parent
                node[above]{$(1,1,1)$}
                node[below]{$f^*$}
                }
                edge from parent
                node[above] {$(1,1)$}
                node[below]  {$f^*$}
            }
            edge from parent
            node[above] {$\widebar{1}$}
            node[below]{${\rm Id}$}
    };
\end{tikzpicture}
}
\end{center}
 With $y_c(t):=  \E\big[ \mathcal{H}(\mathcal{T}_{t,c}) \big]$,
 $c \in \mathcal{C}$,
 the system \eqref{fjhkdsf} reads  
\begin{equation}
\begin{cases}
\nonumber 
\displaystyle
y_{{\rm Id}}(t) 
= y_0 + \int_0^t y_{f^*}(s) ds
\\
\\
\displaystyle
y_{(f^{(k)})^*} (t) = f^{(k)}(y_0) + \int_0^t y_{f^*}(s) y_{(f^{(k+1)})^*}(s) ds, \qquad k \geq 0. 
\end{cases}
\end{equation} 
\section{
Mapping of coding trees to Butcher trees}
\label{s5}
In this section we describe the connection between
coding trees, Butcher series and Butcher trees,
by showing how any Butcher tree
can be recovered by performing a depth first search
on the corresponding coding trees. 
The solution $y(t)$ of the $d$-dimensional ODE system
\eqref{fjklds} is written as the Butcher series 
\begin{align} 
  \nonumber
  y(t) 
  & = 
  y_0 + t f(y_0) + \frac{t^2}{2}  f(y_0)f'(y_0)
  \\
  \nonumber
  & \quad + \frac{t^3}{3!} \big( f''(y_0)f^2(y_0) + f'^2(y_0)f(y_0) \big) 
  \\
  \nonumber
  & \quad + \frac{t^4}{4!} \big( f'''(y_0)f^3(y_0) + f''(y_0)f'(y_0)f^2(y_0)
  + f'(y_0)f'(y_0)f^2(y_0)
  + f'^3(y_0)f(y_0) \big) + \cdots 
  \\
  \label{fhjklsdf}
  & = \quad y_0 + \sum_{\mathcal{B}} \frac{t^{|\mathcal{B}|}}{\nu ( \mathcal{B} )}
  c ({\mathcal{B}})
\end{align} 
where the above summation over Butcher trees ${\mathcal{B}}$
is formal and may not converge.
Here, the order $|\mathcal{B}|$ denotes the
number of vertices of the tree $\mathcal{B}$, 
 $c({\mathcal{B}})$
is a term also depending on the derivatives of $f$ at $y_0$
as in \eqref{aa1}-\eqref{aa3}, 
and $\nu ( \mathcal{B} )$ is a coefficient
which is defined recursively, see Relation~(4) in \cite{butcher2010}. 
\begin{prop}
 Every Butcher tree ${\cal B}$ can be mapped to a unique
 tree ${\rm bin }({\cal B})$, so that 
 the Butcher series \eqref{fhjklsdf} can be rewritten as 
$$
 y_i (t) 
 = y_i(0) + \sum_{\mathcal{B}}
 \E \left[ \mathcal{H}(\mathcal{T}_{t,{\rm Id}_i})
   {\bf 1}_{\big\{
     \mathcal{T}_{t,{\rm Id}_i} \simeq \ \! {\rm bin}(\mathcal{B})
     \big\}
   } 
   \right],
 \quad i=1,\ldots , d.
$$ 
 where the notation $\mathcal{T}_{t,{\rm Id}_i} \simeq \ \! {\rm bin}(\mathcal{B})$
 means that the coding tree $\mathcal{T}_{t,{\rm Id}_i}$ has the tree structure
 ${\rm bin}(\mathcal{B})$. 
\end{prop} 
\begin{Proof}
$(i)$ The one-to-one mapping ${\cal B} \leftrightarrow {\rm bin }({\cal B})$
 is constructed by performing a depth first search on coding trees,
 as illustrated in the following examples 
 which show how any Butcher tree $\mathcal{B}$
 of order $k\geq 1$ 
 can be mapped to a unique tree ${\rm bin}({\cal B})$ 
 generating $k^d$ coding trees $({\cal T}_{i_1,\ldots ,i_k} )_{1 \leq i_1,\ldots ,i_k \leq d}$.
Namely, every Butcher tree
can be recovered by performing a depth first search
on the corresponding coding tree by matching leaves on Butcher trees
to branches of the same color in coding trees.

\medskip

If a node in the initial coding tree branches into two new coding trees ${\cal T}_1$ (above) and ${\cal T}_2$ (below), the Butcher tree for the initial coding tree is obtained by sticking the Butcher tree obtained from ${\cal T}_1$ to the root of the Butcher tree obtained from ${\cal T}_2$.
If a tree is a leaf with code of the form
$c = \big(\partial_{j_1}\cdots \partial_{j_l} f_i \big)^*$ 
 with $i\in \{1,\ldots ,d\}$, 
then the corresponding Butcher tree has a single node containing
$f^{(l)}$. 

\tikzstyle{level 1}=[level distance=4cm, sibling distance=4cm]
\tikzstyle{level 2}=[level distance=4cm, sibling distance=4cm]

\begin{table}[H]
\centering
\scalebox{0.8}{
  \begin{tabular}{||c|C{2.9cm}|C{3cm}|C{10cm}||}
    \hline Order $k$ & $c({\mathcal{B}})$ & Butcher tree $\mathcal{B}$ & Coding trees $({\cal T}_{i_1})_{1\leq i_1\leq d}$ \\
\hline \hline
 1 & $ f (y_0) $ &  
\scalebox{0.8}{
\begin{tikzpicture}[scale=0.7,grow=right, sloped]
\node[circle,draw,black,text=black,thick]{$f$};
\end{tikzpicture}
}&
\scalebox{0.8}{
\begin{tikzpicture}[scale=0.7,grow=right, sloped]
\node[rectangle,draw,black,text=black,thick]{$0$}
    child {
        node[rectangle,draw,black,text=black,thick] {$T_{\widebar{1}}$}
            child {
                node[rectangle,draw,black,text=black,thick] {$t$}
                edge from parent
                node[above] {$(1,1)$}
                node[below]  {$f_{i_1}$}
            }
            edge from parent
            node[above] {$\widebar{1}$}
            node[below]{${\rm Id}_{i_1}$}
    };
\end{tikzpicture}
}\\
\hline\hline
\end{tabular}
}
\end{table}

\vspace{-0.6cm}

\begin{table}[H]
\centering
\scalebox{0.8}{
  \begin{tabular}{||c|C{2.9cm}|C{3cm}|C{10cm}||}
    \hline Order $k$ & $c({\mathcal{B}})$ & Butcher tree $\mathcal{B}$ & Coding trees $({\cal T}_{i_1,i_2})_{1\leq i_1,i_2\leq d}$ \\
\hline \hline
 2 & $ff'(y_0)$ &
\scalebox{0.8}{
\begin{tikzpicture}[scale=0.4,grow=up, sloped]
\node[circle,draw,black,text=black,thick]{$f'$}
    child[draw=blue]{
        node[circle,draw,blue,text=blue,thick] {$f$}
            edge from parent
    };
\end{tikzpicture}}&
\scalebox{0.7}{
\begin{tikzpicture}[scale=0.7,grow=right, sloped]
\node[rectangle,draw,black,text=black,thick]{$0$}
    child {
        node[rectangle,draw,black,text=black,thick] {$T_{\widebar{1}}$}
            child {
                node[rectangle,draw,black,text=black,thick] {$T_{(1,1)}$}
                child{
                node[rectangle,draw,black,text=black,thick]{$t$}
                edge from parent
                node[above]{$(1,1,2)$}
                node[below]{$\partial_{i_2}f_{i_1}$}
                }
                child[draw=blue]{
                node[rectangle,black,draw,thick]{$t$}
                edge from parent
                node[above]{$(1,1,1)$}
                node[below]{$\textcolor{blue}{f}_{i_2}$}
                }
                edge from parent
                node[above] {$(1,1)$}
                node[below]  {$f_{i_1}$}
            }
            edge from parent
            node[above] {$\widebar{1}$}
            node[below]{${\rm Id}_{i_1}$}
    };
\end{tikzpicture}
}\\
\hline\hline
\end{tabular}
}
\end{table}

\vspace{-0.6cm}

\begin{table}[H]
\centering
\scalebox{0.8}{
  \begin{tabular}{||c|C{2.9cm}|C{3cm}|C{10cm}||}
    \hline
Order $k$ & $c({\mathcal{B}})$ & Butcher tree $\mathcal{B}$ & Coding trees $({\cal T}_{i_1,i_2,i_3})_{1\leq i_1,i_2,i_3\leq d}$ \\
\hline \hline
3 & $f''f^2(y_0)$ &
\scalebox{0.8}{
\begin{tikzpicture}[scale=0.4,grow=up, sloped]
\node[circle,draw,black,text=black,thick]{$f''$}
    child[draw=purple] {
        node[circle,draw,purple,text=black,thick] {$\textcolor{purple}{f}$}
            edge from parent
    }
    child[draw=blue] {
        node[circle,draw,blue,text=black,thick] {$\textcolor{blue}{f}$}
            edge from parent
    }
    ;
\end{tikzpicture}}&
\scalebox{0.7}{
\begin{tikzpicture}[scale=0.7,grow=right, sloped]
\node[rectangle,draw,black,text=black,thick]{$0$}
    child {
        node[rectangle,draw,black,text=black,thick] {$T_{\widebar{1}}$}
            child {
                node[rectangle,draw,black,text=black,thick] {$T_{(1,1)}$}
                child{
                node[rectangle,draw,black,text=black,thick]{$T_{(1,1,2)}$}
                child{
                    node[rectangle,draw,black,text=black,thick]{$t$}
                    edge from parent
                    node[above]{$(1,1,2,2)$}
                    node[below]{$\partial_{i_3}\partial_{i_2}f_{i_1}$}
                }
                child[draw=blue]{
                    node[rectangle,draw,black,text=black,thick]{$t$}
                    edge from parent
                    node[above]{$(1,1,2,1)$}
                    node[below]{$\textcolor{blue}{f_{i_3}}$}
                }
                edge from parent
                node[above]{$(1,1,2)$}
                node[below]{$\partial_{i_2}f_{i_1}$}
                }
                child[draw=purple]{
                node[rectangle,black,draw,thick]{$t$}
                edge from parent
                node[above]{$(1,1,1)$}
                node[below]{$\textcolor{purple}{f_{i_2}}$}
                }
                edge from parent
                node[above] {$(1,1)$}
                node[below]  {$f_{i_1}$}
            }
            edge from parent
            node[above] {$\widebar{1}$}
            node[below]{${\rm Id}_{i_1}$}
    };
\end{tikzpicture}
}
\\
\hline\hline
\end{tabular}
}
\end{table}

\vspace{-0.6cm}

\begin{table}[H]
\centering
\scalebox{0.8}{
\begin{tabular}{||C{1.4cm}|C{2.9cm}|C{3cm}|C{10cm}||}
\hline
\hline
3 & $f'^2f(y_0)$ &
\scalebox{0.8}{
\begin{tikzpicture}[scale=0.3,grow=up, sloped]
\node[circle,draw,black,text=black,thick]{$f'$}
    child[draw=blue] {
        node[circle,draw,blue,text=blue,thick] {$f'$}
        child[draw=purple] {
            node[circle,draw,purple,text=purple,thick] {$f$}
                edge from parent
        edge from parent
    }
    }
    ;
\end{tikzpicture}
}
&
\scalebox{0.65}{
\begin{tikzpicture}[scale=0.8,grow=right, sloped]
\node[rectangle,draw,black,text=black,thick]{$0$}
    child {
        node[rectangle,draw,black,text=black,thick] {$T_{\widebar{1}}$}
            child {
                node[rectangle,draw,black,text=black,thick] {$T_{(1,1)}$}
                child{
                node[rectangle,draw,black,text=black,thick]{$t$}
                edge from parent
                node[above]{$(1,1,2)$}
                node[below]{$\partial_{i_2}f_{i_1}$}
                }
                child{
                node[rectangle,draw,thick]{$T_{(1,1,1)}$}
                    child[draw=blue]{
                    node[rectangle,draw,black,text=black,thick]{$t$}
                    edge from parent
                    node[above]{$(1,1,1,2)$}
                    node[below]{$\textcolor{blue}{\partial_{i_3}f_{i_2}}$}
                    }
                    child[draw=purple]{
                    node[rectangle,draw,black,text=black,thick]{$t$}
                    edge from parent
                    node[above]{$(1,1,1,1)$}
                    node[below]{$\textcolor{purple}{f_{i_3}}$}
                    }
                edge from parent
                node[above]{$(1,1,1)$}
                node[below]{$f_{i_2}$}
                }
                edge from parent
                node[above] {$(1,1)$}
                node[below]  {$f_{i_1}$}
            }
            edge from parent
            node[above] {$\widebar{1}$}
            node[below]{${\rm Id}_{i_1}$}
    };
\end{tikzpicture}
}
\\
\hline\hline
\end{tabular}
}
\end{table}

\vspace{-0.6cm}

\begin{table}[H]
\centering
\scalebox{0.8}{
           {
            \begin{tabular}{||c|C{2.9cm}|C{3cm}|C{10cm}||}
              \hline
Order $k$ & $c({\mathcal{B}})$ & Butcher tree $\mathcal{B}$ & Coding trees $({\cal T}_{i_1,i_2,i_3,i_4})_{1\leq i_1,i_2,i_3,i_4\leq d}$ \\ 
\hline \hline
 4 & $f'''f^3(y_0)$ &
\scalebox{0.8}{
\begin{tikzpicture}[scale=0.3,grow=up, sloped]
\node[circle,draw,black,text=black,thick]{$f'''$}
    child [draw=cyan]{
        node[circle,draw,cyan,text=cyan,thick] {$f$}
            edge from parent
    }
    child[draw=purple] {
        node[circle,draw,purple,text=purple,thick] {$f$}
            edge from parent
    }
    child[draw=blue] {
        node[circle,draw,blue,text=blue,thick] {$f$}
            edge from parent
    }
    ;
\end{tikzpicture}
          }
          &
\scalebox{0.6}{
\begin{tikzpicture}[scale=0.7,grow=right, sloped]
\node[rectangle,draw,black,text=black,thick]{$0$}
    child {
        node[rectangle,draw,black,text=black,thick] {$T_{\widebar{1}}$}
            child {
                node[rectangle,draw,black,text=black,thick] {$T_{(1,1)}$}
                child{
                node[rectangle,draw,black,text=black,thick]{$T_{(1,1,2)}$}
                    child{
                        node[rectangle,draw,black,text=black,thick]{$T_{(1,1,2,2)}$}
                        child{
                            node[rectangle,draw,black,text=black,thick]{$t$}
                        edge from parent
                        node[above]{$(1,1,2,2,1)$}
                        node[below]{$\partial_{i_4}\partial_{i_3}\partial_{i_2}f_{i_1}$}
                        }
                        child[draw=blue]{
                            node[rectangle,draw,black,text=black,thick]{$t$}
                        edge from parent
                        node[above]{$(1,1,2,2,2)$}
                        node[below]{$\textcolor{blue}{f_{i_4}}$}
                        }
                    edge from parent
                    node[above]{$(1,1,2,2)$}
                    node[below]{$\partial_{i_3}\partial_{i_2}f_{i_1}$}
                    }
                    child[draw=purple]{
                        node[rectangle,draw,black,text=black,thick]{$t$}
                        edge from parent
                        node[above]{$(1,1,2,1)$}
                        node[below]{$\textcolor{purple}{f_{i_3}}$}
                    }
                edge from parent
                node[above]{$(1,1,2)$}
                node[below]{$\partial_{i_2}f_{i_1}$}
                }
                child[draw=cyan]{
                node[rectangle,draw,thick]{$t$}
                edge from parent
                node[above]{$(1,1,1)$}
                node[below]{$\textcolor{cyan}{f_{i_2}}$}
                }
                edge from parent
                node[above] {$(1,1)$}
                node[below]  {$f_{i_1}$}
            }
            edge from parent
            node[above] {$\widebar{1}$}
            node[below]{${\rm Id}_{i_1}$}
    };
\end{tikzpicture}
       }
          \\
\hline\hline
\end{tabular}
          }
          }
\end{table}

\vspace{-0.6cm}

\begin{table}[H]
\centering
\scalebox{0.8}{
\begin{tabular}{||C{1.4cm}|C{2.9cm}|C{3cm}|C{10cm}||}
\hline \hline
4 & $f''f'f^2(y_0)$ &
\scalebox{0.8}{
\begin{tikzpicture}[scale=0.3,grow=up, sloped]
\node[circle,draw,black,text=black,thick]{$f''$}
    child[draw=purple] {
        node[circle,draw,purple,text=purple,thick] {$f'$}
            child[draw=cyan]{
                    node[circle,draw,cyan,text=cyan,thick] {$f$}
            edge from parent
            }
            edge from parent
    }
    child[draw=blue]{
        node[circle,draw,blue,text=blue,thick] {$f$}
            edge from parent
    }
    ;
\end{tikzpicture}
}
&
\scalebox{0.6}{
\begin{tikzpicture}[scale=0.8,grow=right, sloped]
\node[rectangle,draw,black,text=black,thick]{$0$}
    child {
        node[rectangle,draw,black,text=black,thick] {$T_{\widebar{1}}$}
            child {
                node[rectangle,draw,black,text=black,thick] {$T_{(1,1)}$}
                child{
                node[rectangle,draw,black,text=black,thick,yshift=-0.4cm]{$T_{(1,1,2)}$}
                    child{
                        node[rectangle,draw,black,text=black,thick]{$t$}
                    edge from parent
                    node[above]{$(1,1,2,2)$}
                    node[below]{$\partial_{i_3}\partial_{i_2}f_{i_1}$}
                    }
                    child[draw=blue]{
                        node[rectangle,draw,black,text=black,thick]{$t$}
                        edge from parent
                        node[above]{$(1,1,2,1)$}
                        node[below]{$\textcolor{blue}{f_{i_3}}$}
                    }
                edge from parent
                node[above]{$(1,1,2)$}
                node[below]{$\partial_{i_2}f_{i_1}$}
                }
                child{
                node[rectangle,draw,thick,yshift=0.4cm]{$T_{(1,1,1)}$}
                    child[draw=purple]{
                    node[rectangle,draw,black,text=black,thick]{$t$}
                    edge from parent
                    node[above]{$(1,1,1,2)$}
                    node[below]{$\textcolor{purple}{\partial_{i_4}f_{i_2}}$}
                    }
                    child[draw=cyan]{
                    node[rectangle,draw,black,text=black,thick]{$t$}
                    edge from parent
                    node[above]{$(1,1,1,1)$}
                    node[below]{$\textcolor{cyan}{f_{i_4}}$}
                    }
                edge from parent
                node[above]{$(1,1,1)$}
                node[below]{$f_{i_2}$}
                }
                edge from parent
                node[above] {$(1,1)$}
                node[below]  {$f_{i_1}$}
            }
            edge from parent
            node[above] {$\widebar{1}$}
            node[below]{${\rm Id}_{i_1}$}
    };
\end{tikzpicture}
       }
          \\
\hline\hline
\end{tabular}
}
\end{table}

\vspace{-0.6cm}

\begin{table}[H]
\centering
\scalebox{0.8}{
\begin{tabular}{||C{1.4cm}|C{2.9cm}|C{3cm}|C{10cm}||}
\hline \hline
4 & $f'f''f^2(y_0)$ &
\scalebox{0.8}{
\begin{tikzpicture}[scale=0.3,grow=up, sloped]
\node[circle,draw,black,text=black,thick]{$f'$}
    child[draw=blue] {
        node[circle,draw,blue,text=blue,thick] {$f''$}
            child[draw=cyan]{
                    node[circle,draw,cyan,text=cyan,thick] {$f$}
            edge from parent
            }
            child[draw=purple]{
                    node[circle,draw,purple,text=purple,thick] {$f$}
            edge from parent
            }
            edge from parent
    }
    ;
\end{tikzpicture}
}
&
\scalebox{0.6}{
\begin{tikzpicture}[scale=0.8,grow=right, sloped]
\node[rectangle,draw,black,text=black,thick]{$0$}
    child {
        node[rectangle,draw,black,text=black,thick] {$T_{\widebar{1}}$}
            child {
                node[rectangle,draw,black,text=black,thick] {$T_{(1,1)}$}
                child{
                    node[rectangle,draw,black,text=black,thick]{$t$}
                edge from parent
                node[above]{$(1,1,2)$}
                node[below]{$\partial_{i_2} f_{i_1}$}
                }
                child{
                node[rectangle,draw,thick]{$T_{(1,1,1)}$}
                    child{
                    node[rectangle,draw,black,text=black,thick]{$T_{(1,1,1,2)}$}
                    child[draw=blue]{
                        node[rectangle,draw,black,text=black,thick]{$t$}
                    edge from parent
                    node[above]{$(1,1,1,2,2)$}
                    node[below]{$\textcolor{blue}{\partial_{i_4}\partial_{i_3} f_{i_2}}$}                    
                    }
                    child[draw=purple]{
                        node[rectangle,draw,black,text=black,thick]{$t$}
                    edge from parent
                    node[above]{$(1,1,1,2,1)$}
                    node[below]{$\textcolor{purple}{f_{i_4}}$}                    
                    }
                    edge from parent
                    node[above]{$(1,1,1,2)$}
                    node[below]{$\partial_{i_3}f_{i_2}$}
                    }
                    child[draw=cyan]{
                    node[rectangle,draw,black,text=black,thick]{$t$}
                    edge from parent
                    node[above]{$(1,1,1,1)$}
                    node[below]{$\textcolor{cyan}{f_{i_3}}$}
                    }
                edge from parent
                node[above]{$(1,1,1)$}
                node[below]{$f_{i_2}$}
                }
                edge from parent
                node[above] {$(1,1)$}
                node[below]  {$f_{i_1}$}
            }
            edge from parent
            node[above] {$\widebar{1}$}
            node[below]{${\rm Id}_{i_1}$}
    };
\end{tikzpicture}
       }
          \\
\hline\hline
\end{tabular}
}
\end{table}

\vspace{-0.6cm}

\begin{table}[H]
\centering
\scalebox{0.8}{
\begin{tabular}{||C{1.4cm}|C{2.9cm}|C{3cm}|C{10cm}||}
\hline \hline
 4 & $f'^3f(y_0)$ &
\scalebox{0.8}{
\begin{tikzpicture}[scale=0.3,grow=up, sloped]
\node[circle,draw,black,text=black,thick]{$f'$}
    child {
        node[circle,draw,blue,text=blue,thick] {$f'$}
            child{
                    node[circle,draw,purple,text=purple,thick] {$f'$}
                        child{
                            node[circle,draw,cyan,text=cyan,thick] {$f$}
                        edge from parent
                        }
            edge from parent
            }
            edge from parent
    }
    ;
\end{tikzpicture}
}
&
\scalebox{0.6}{
\begin{tikzpicture}[scale=0.8,grow=right, sloped]
\node[rectangle,draw,black,text=black,thick]{$0$}
    child {
        node[rectangle,draw,black,text=black,thick] {$T_{\widebar{1}}$}
            child {
                node[rectangle,draw,black,text=black,thick] {$T_{(1,1)}$}
                child{
                    node[rectangle,draw,black,text=black,thick]{$t$}
                edge from parent
                node[above]{$(1,1,2)$}
                node[below]{$\partial_{i_2}f_{i_1}$}
                }
                child{
                node[rectangle,draw,thick]{$T_{(1,1,1)}$}
                    child[draw=blue]{
                    node[rectangle,draw,black,text=black,thick]{$t$}
                    edge from parent
                    node[above]{$(1,1,1,2)$}
                    node[below]{$\textcolor{blue}{\partial_{i_3}f_{i_2}}$}
                    }
                    child{
                    node[rectangle,draw,black,text=black,thick]{$T_{(1,1,1,1)}$}
                        child[draw=purple]{
                            node[rectangle,draw,black,text=black,thick]{$t$}
                        edge from parent
                        node[above]{$(1,1,1,1,2)$}
                        node[below]{$\textcolor{purple}{\partial_{i_4} f_{i_3}}$}
                        }
                        child[draw=cyan]{
                            node[rectangle,draw,black,text=black,thick]{$t$}
                        edge from parent
                        node[above]{$(1,1,1,1,1)$}
                        node[below]{$\textcolor{cyan}{f_{i_4}}$}
                        }
                    edge from parent
                    node[above]{$(1,1,1,1)$}
                    node[below]{$f_{i_3}$}
                    }
                edge from parent
                node[above]{$(1,1,1)$}
                node[below]{$f_{i_2}$}
                }
                edge from parent
                node[above] {$(1,1)$}
                node[below]  {$f_{i_1}$}
            }
            edge from parent
            node[above] {$\widebar{1}$}
            node[below]{${\rm Id}_{i_1}$}
    };
\end{tikzpicture}
}
\\
\hline\hline
\end{tabular}
}
\end{table}
        
\vspace{-0.4cm}

\noindent
$(ii)$
 When ${\mathcal{B}}$ is a Butcher tree of order $k=|{\cal B}| \geq 1$,
 based on the construction of $\nu ({\cal B}) = \sigma ({\cal B})
 \gamma ( {\cal B})$ on page~155 of \cite{butcher2010}, 
 the $i_1$-$th$ component of the term
 $t^{|\mathcal{B}|} c({\mathcal{B}}) / \nu ( \mathcal{B} )$ in the expansion of $y_{i_1}(t)$ 
 is a sum of $f$ and its partial derivatives at $y_0$ as
 in \eqref{aa1}-\eqref{aa3}.
 This sum can be interpreted as 
 a sum of $\mathcal{H}( {\cal T}_{i_1,\ldots ,i_k} )$
 over $1 \leq i_2,\ldots ,i_k \leq d$,
 weighted by the probability coefficients 
$$
\displaystyle \prod_{\widebar{k} \in \mathcal{K}^{\circ}} \big(
q_{c_{\widebar{k}}}(I_{c_{\widebar{k}}})\rho(T_{\widebar{k}}-T_{\widebar{k}{\text{--}}})
\big)
\prod_{\widebar{k} \in \mathcal{K}^{\partial}}
\big( \widebar{F}(t-T_{\widebar{k}{\text{--}}})\big)
$$
and integrated over the branching times $T_{\widebar{k}} \in [0,t]$,
which yields the identity 
$$
\frac{t^{|\mathcal{B}|}}{\nu ( \mathcal{B} )}
c ({\mathcal{B}})
=
\left( \E \left[ \mathcal{H}(\mathcal{T}_{t,{\rm Id}_{i_1}})
   {\bf 1}_{\big\{
     \mathcal{T}_{t,{\rm Id}_{i_1}} \simeq \ \! {\rm bin}(\mathcal{B})
     \big\}
   } 
   \right]
\right)_{i_1=1,\ldots , d}.
$$ 
\end{Proof} 
\section{Numerical application} 
\label{s6}
In this section we consider various ODE examples
with different choices of probability density functions
$\rho (t)$ satisfying \eqref{fjkld}, or under \eqref{p1.1}. 
 The errors observed in Figures~\ref{f6}, \ref{f3}, \ref{f4},
 only occur at a 
 time threshold after which the estimator
 $\mathcal{H}(\mathcal{T}_{t,c})$ is no longer integrable and
 the estimates are unreliable.
 The following graphs are plotted with one million Monte Carlo samples.
 
\begin{enumerate}[i)]
\item Taking $f(y):=y^2$, we start with the quadratic ODE
  \begin{equation}
    \label{e1} 
y'(t) = y^2(t), \qquad y(0) = y_0 = 1, 
\end{equation} 
 with solution
$$
y(t) = \frac{y_0}{1-y_0 t}, \qquad t\in [0,1/y_0). 
$$
   In the framework of \eqref{Eint2} we have 
 $$
 \mathcal{C}= \big\{{\rm Id}, \ \big( f^{(k)} \big)^*, \ k\geq 0\big\}
 = \big\{{\rm Id}, \ (x\mapsto 0)^*, \ (x\mapsto 2)^*, \ (x\mapsto 2x)^*, \ (x\mapsto x^2)^* \big\},
 $$ 
 hence we have 
 $c(y)(0) \leq K:=\max ( 1 , 2y_0, y_0^2) = 2$ for all $c\in \mathcal{C}$,
 and in agreement with \eqref{p1.1}   
 the representation formula \eqref{fjhkldsf} of Theorem~\ref{t1} holds
 for $t \in [0,1/2)$, see Figure~\ref{f1}.
In this example and in the next two examples we take $\rho$ to be the exponential 
probability density function $\rho (t)=e^{-t}$, $t\geq 0$.

\begin{figure}[H]
\centering
\includegraphics[width=0.6\textwidth]{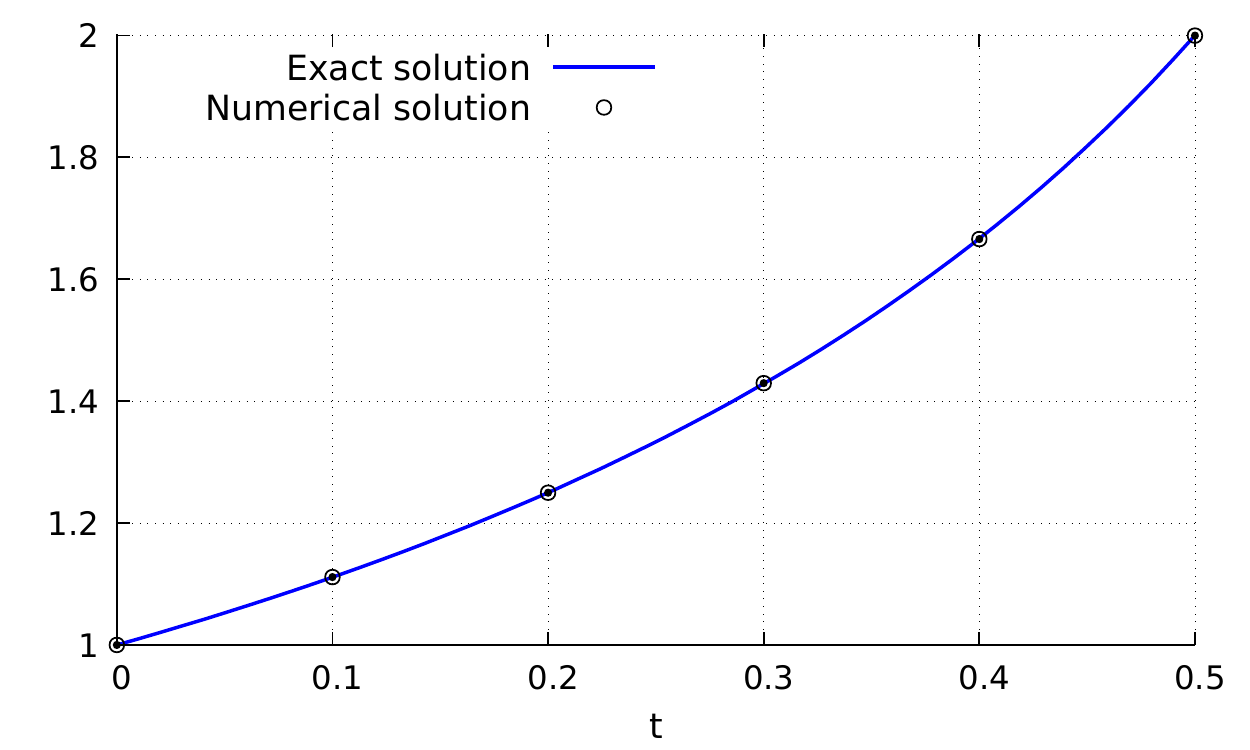}
\vskip-0.2cm
\caption{Numerical solution of \eqref{e1} with $y_0=1$.}
\label{f1}
\end{figure}

\vspace{-0.6cm}

\item Next, we take $f(y) := \cos(y)$ and consider the equation 
   \begin{equation}
     \label{e2} 
y'(t) = \cos(y(t)), \quad y(0) = y_0, 
\end{equation} 
 with solution
$$
 y(t) = 2 \tan^{-1}\left(\tanh\left(\frac{t+2\tanh^{-1}(\tan( y_0 / 2 ))}{2}\right)\right),
  \qquad t\in \real_+. 
$$
in the framework of \eqref{Eint2}.
When $y_0=1$ we have $K = \sup_{k\in\mathbb{N}} f^{(k)}(1) = 1$, 
 and in agreement with \eqref{p1.1}    
 the representation formula \eqref{fjhkldsf} of
 Theorem~\ref{t1} holds
 for $t \in [0,1)$, see Figure~\ref{f2}.

\begin{figure}[H]
\centering
\includegraphics[width=0.6\textwidth]{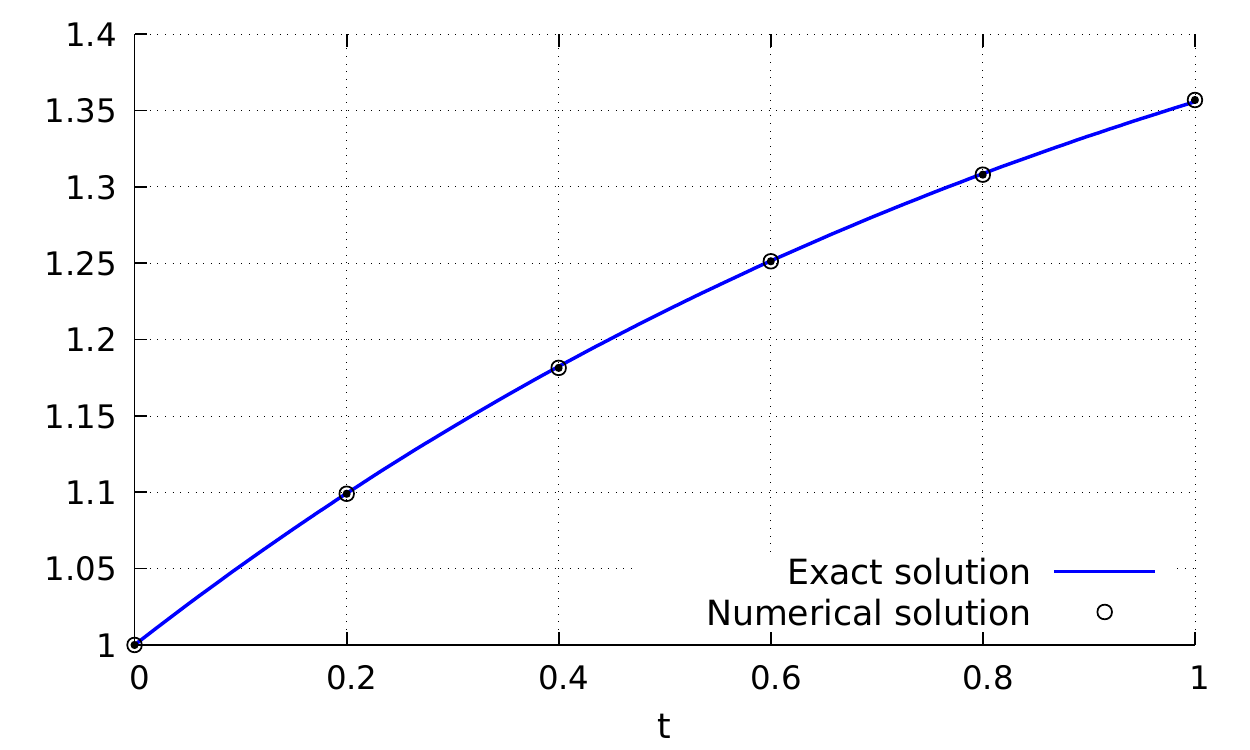}
\vskip-0.2cm
\caption{Numerical solution of \eqref{e2}.}
\label{f2}
\end{figure}

\vspace{-0.6cm}

\item
  \label{ex3}
  Taking $f(t,y) := (y+t)/(y-t)$ we find
  Equation~(201a) in \cite{butcherbk}, i.e. 
  \begin{equation}
    \label{e4.0} 
y'(t) = \frac{y(t)+t}{y(t)-t},  \qquad y(0) = 1, 
\end{equation} 
 with solution 
$$
y(t) = t+\sqrt{1+2t^2}. 
$$
In this case, the time interval of validity
may not be determined explicitly 
because $\sup_{k,l\geq 0} |\partial_0^{k}\partial_1^{l}f(0,1/2)|=\infty$
and the finiteness of this supremum is only a sufficient condition
for \eqref{gfdg} to hold in Proposition~\ref{p1}. 
Figure~\ref{f6} shows the convergence of the Monte Carlo algorithm until $t=0.25$. 

\begin{figure}[H]
\centering
\includegraphics[width=0.6\textwidth]{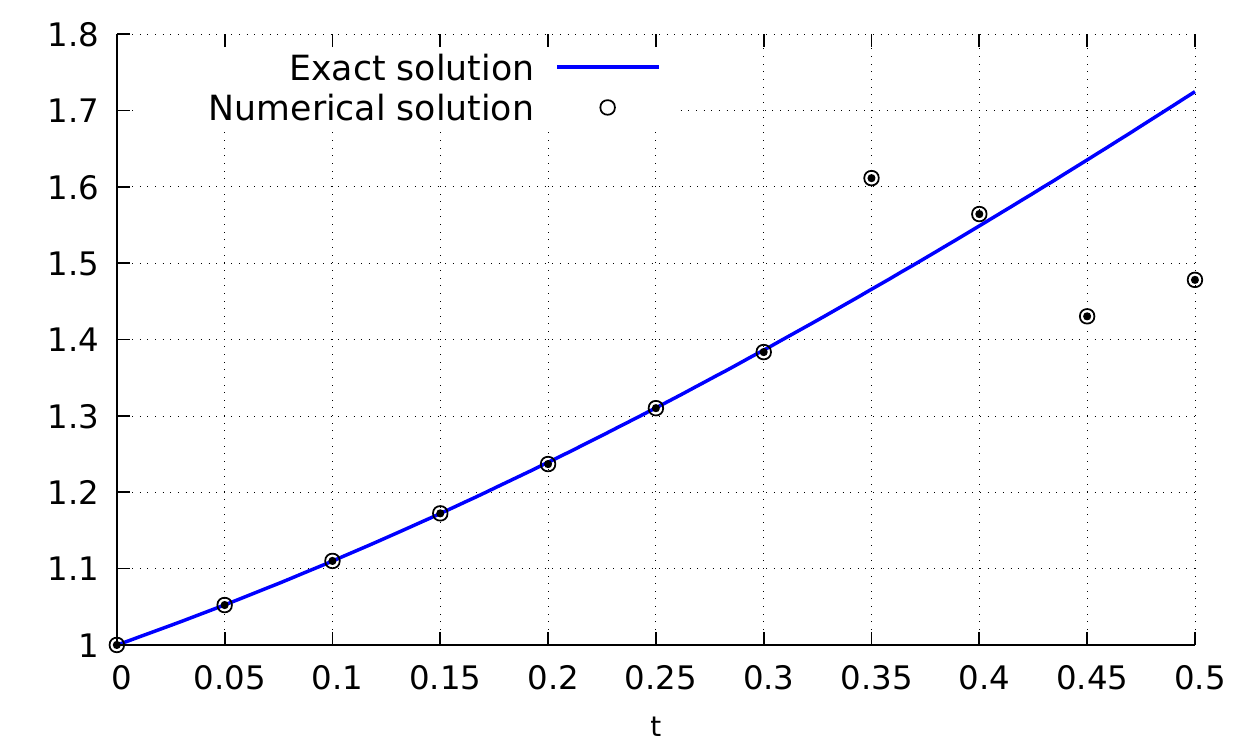}
\vskip-0.2cm
\caption{Numerical solution of \eqref{e4.0}.}
\label{f6}
\end{figure}

\vspace{-0.6cm}

\item Taking $f(t,y) := (y-t)/(y+t)$ yields 
   Equation~(316e) in \cite{butcherbk}, i.e.  
\begin{equation} 
  \label{e3}
  y'(t) = \frac{y(t)-t}{y(t)+t}, \qquad y(0) = 1, 
\end{equation} 
whose solution is given in parametric form as
$(t(u) , y(u) ) = ( u \sin \log (u), u \cos \log (u))$.
 As in Example~\ref{ex3}) above, the time interval of validity
 may not be determined explicitly, see Figure~\ref{f3}. 
In this example and in the next one, we take $\rho$ to be the gamma
probability density function $\rho (t)= t^{-1/2}e^{-t}/\Gamma (1/2)$, $t > 0$.

\begin{figure}[H]
\centering
\includegraphics[width=0.6\textwidth]{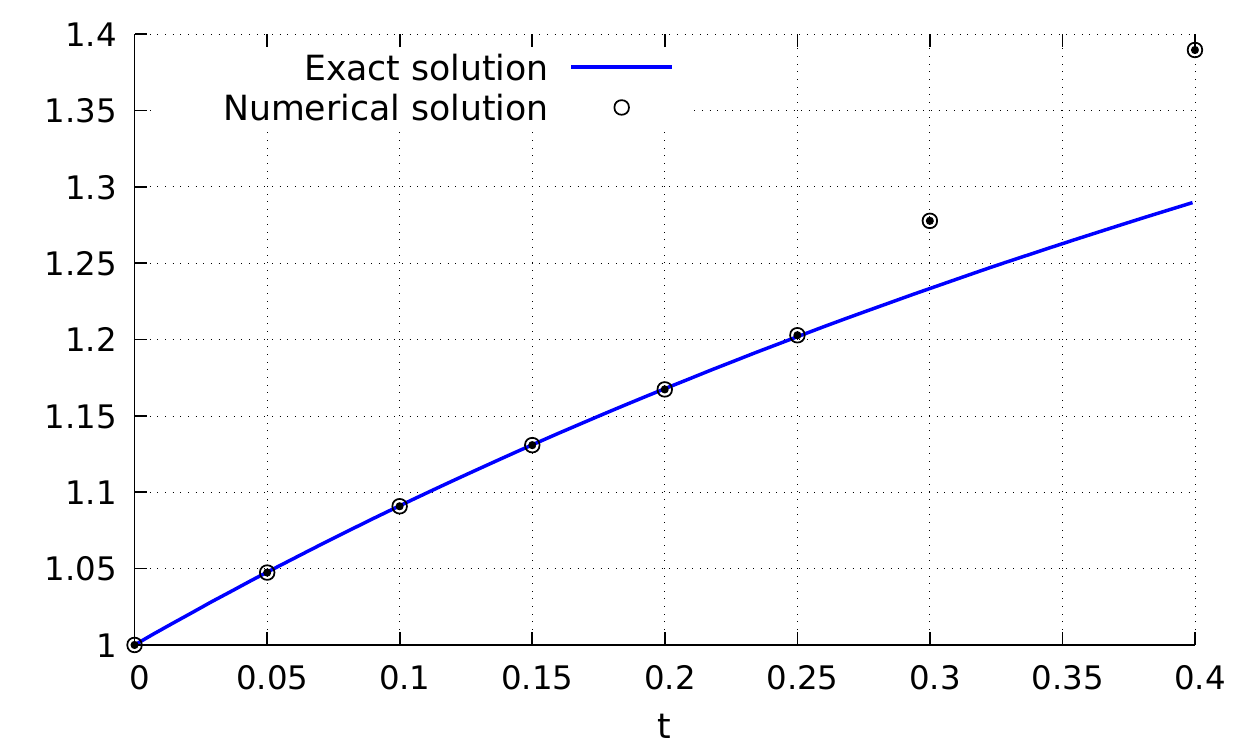}
\vskip-0.2cm
\caption{Numerical solution of \eqref{e3}.}
\label{f3}
\end{figure}

\vspace{-0.6cm}

\item
  Taking $f(t,y) := yt + y^2$ we find
  Equation~(223a) in \cite{butcherbk}, i.e. 
  \begin{equation}
    \label{e4} 
y'(t) = t y(t) + y^2(t),  \qquad y(0) = 1/2, 
\end{equation} 
 with solution 
$$
y(t) = \frac{e^{t^2/2}}{2-\int_0^t e^{s^2/2} ds}. 
$$
 In this case we have $K = \sup_{k,l\geq 0} |\partial_0^{k}\partial_1^{l}f(0,1/2)| \leq 2$, 
 and in agreement with \eqref{p1.1}     
 the representation formula \eqref{fjhkldsf} of
 Theorem~\ref{t1} is valid on the time interval $[0,0.5) \subset [0,1/K]$,
 see Figure~\ref{f4}.

\begin{figure}[H]
\centering
\includegraphics[width=0.6\textwidth]{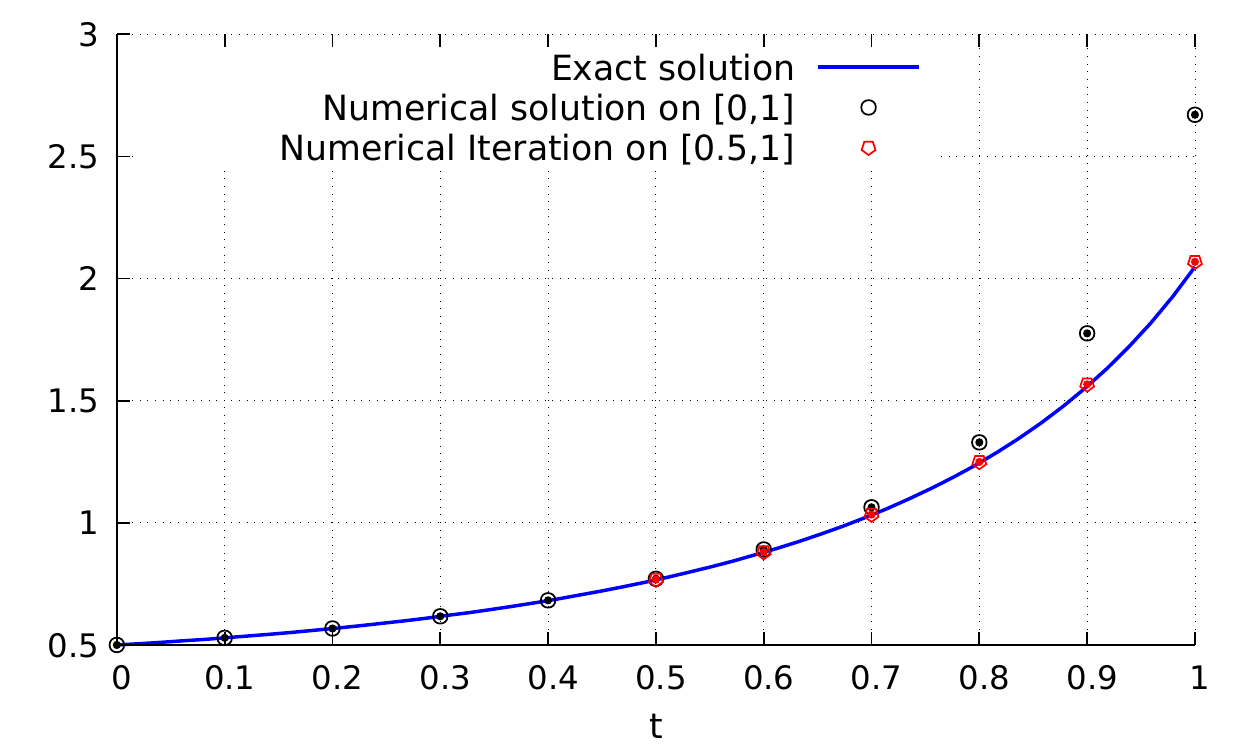}
\vskip-0.2cm
\caption{Numerical solution of \eqref{e4}.}
\label{f4}
\end{figure}

\vspace{-0.6cm}

In addition, after running the algorithm on the time interval $[0,0.5]$
we may reuse the numerical evaluation at time $t=0.5$ as a new initial
condition and iterate the algorithm on the time interval $[0.5, 1]$ with more
precise estimates, as shown in red in Figure~\ref{f4}. We refer to this procedure as
``patching''. 

\item
  Consider the ODE system (316f) page~177 in \cite{butcherbk}, i.e. 
\begin{equation}
\label{e5-ab} 
\left\{
\begin{array}{l}
  y'_1(t) = \displaystyle \frac{y_1(t)+y_2(t)}{\sqrt{y_1^2(t)+y_2^2(t)}},  \qquad y_1(1) = 0, 
  \medskip
  \\
y'_2(t) = \displaystyle \frac{y_2(t)-y_1(t)}{\sqrt{y_1^2(t)+y_2^2(t)}},  \qquad y_2(1) = 1, 
\smallskip
\end{array}
\right.
\end{equation} 
 with solution 
 $$
  y_1(t) = t \sin \log t, 
\qquad 
  y_2(t) = t \cos \log t, \qquad t\geq 1. 
$$
 The graphs of Figure~\ref{fig1}
 are obtained using the Python code provided in appendix,
 where we patch the algorithm $6$ times over the time interval $[1,4]$ 
 and take $\rho$ to be the exponential probability density function
 $\rho (t)= e^{-t}$, $t > 0$.

\begin{figure}[H]
\centering
\begin{subfigure}{.49\textwidth}
\centering
\includegraphics[width=\textwidth]{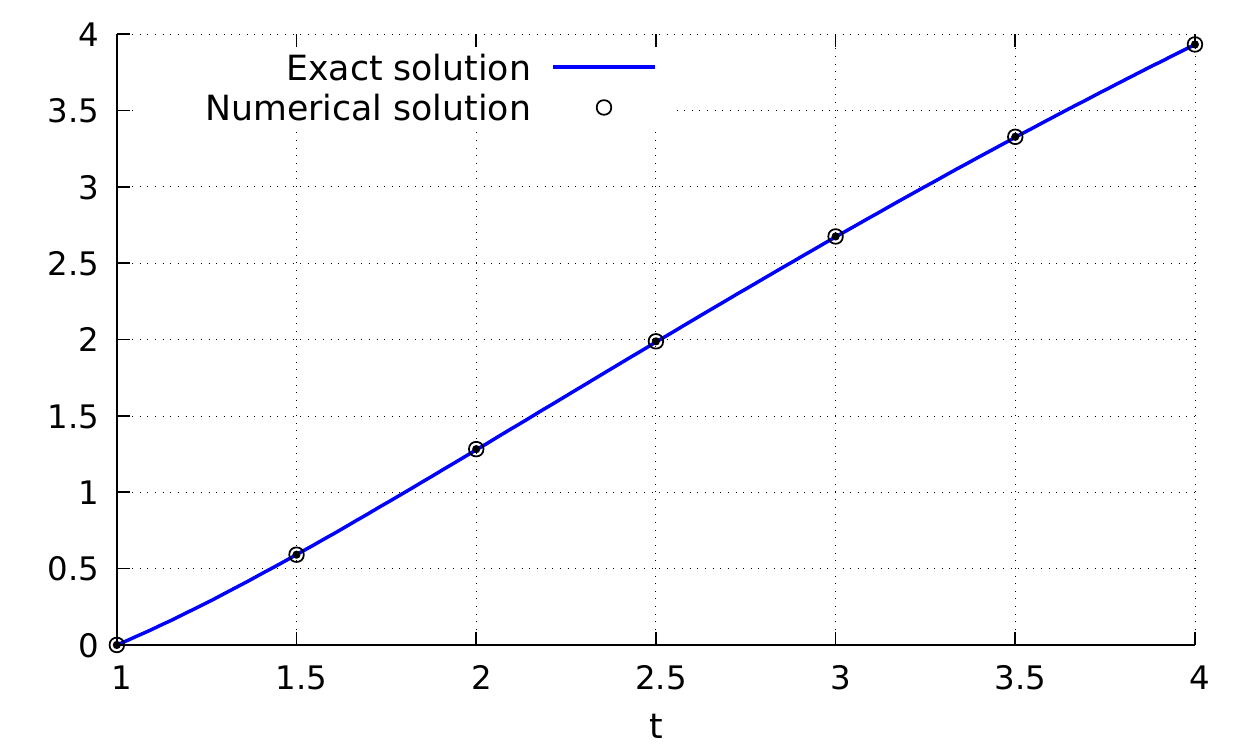}
\vskip-0.1cm
\caption{Graph of $y_1(t)$.}
\end{subfigure}
\hskip-0.2cm
\begin{subfigure}{.49\textwidth}
\centering
\includegraphics[width=\textwidth]{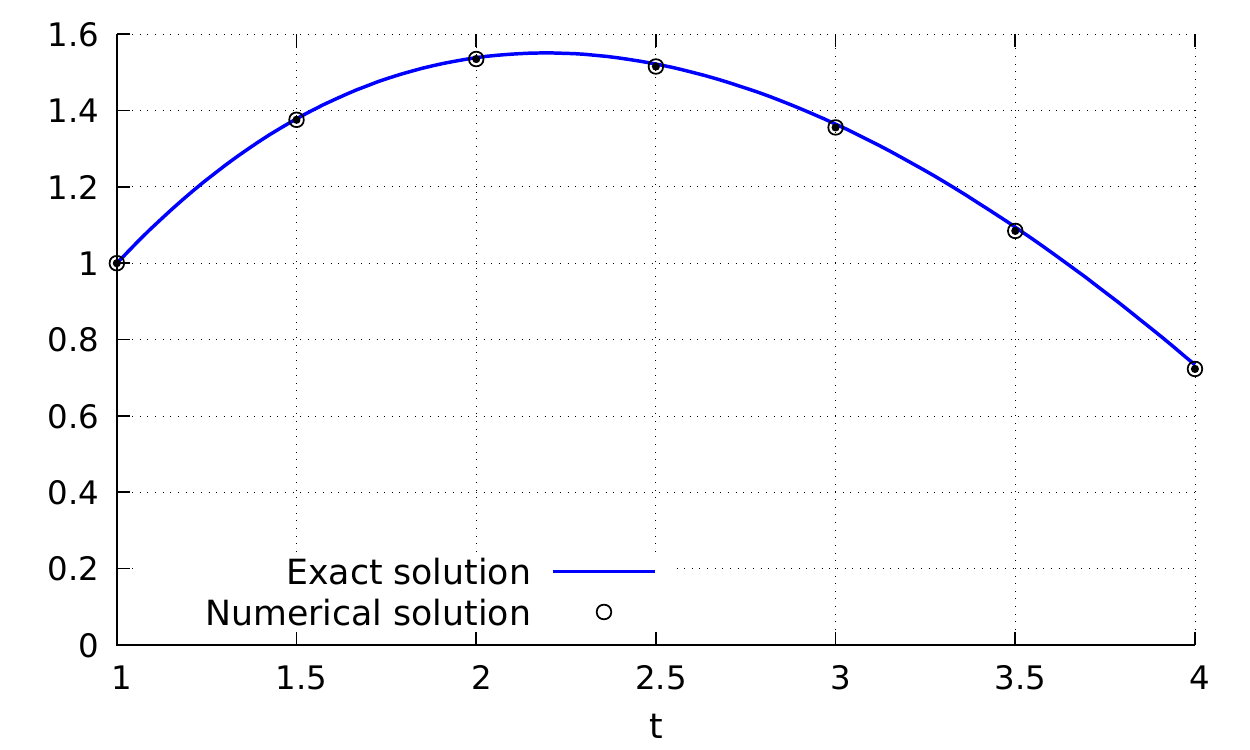}
\vskip-0.1cm
\caption{Graph of $y_2(t)$.}
\end{subfigure}
\caption{Numerical solution of \eqref{e5-ab}.}
\label{fig1}
\end{figure}

\vskip-0.3cm

\end{enumerate}

\appendix

\section{Computer codes} 

The following Maple and Mathematica codes implement the algorithm of Theorem~\ref{t1} for one-dimensional non-autonomous ODEs using the exponential distribution $\rho (t) = \re^{-t}$, $t\geq 0$. 

\medskip

\begin{lstlisting}[language=Maple, caption={1}, title={Maple code 
      .} 
  ]
codetofunction := proc(f, c, t0, y0) if nops(c) = 0 then return y0; end if; if c = [0, 0] then return f(t0, y0); else return eval(eval(diff(f(t, y), t $ c[1], y $ c[2]), t = t0), y = y0); end if; end proc;
mcsample := proc(f, t, t0, y0, c, h) local A, tau; tau := random[exponential[1]](1); if t - t0 < tau then return h*codetofunction(f, c, t0, y0)/exp(-t + t0); else if nops(c) = 0 then return mcsample(f, t - tau, t0, y0, [0, 0], h/exp(-tau)); else if random[uniform](1) < 0.5 then return mcsample(f, t - tau, t0, y0, [c[1] + 1, c[2]], 2*h/exp(-tau)); else A := mcsample(f, t - tau, t0, y0, [0, 0], 1); return mcsample(f, t - tau, t0, y0, [c[1], c[2] + 1], 2*A*h/exp(-tau)); end if; end if; end if; end proc;
solution := proc(f, t, t0, y0, n) local i, temp; temp := 0; for i to n do temp := temp + mcsample(f, t, t0, y0, [], 1); end do; return evalf(temp/n); end proc;
f := (t, y) -> y^2;solution(f, 0.5, 0, 1, 10000);
\end{lstlisting}

\begin{lstlisting}[language=Mathematica,caption={2}, title={Mathematica code.} 
  ]
codetofunction[f_, c__, t0_, y0_] := (If [c == {}, Return [y0], 
   Return[D[D[f[t, y], {t, c[[1]]}], {y, c[[2]]}] /. {t -> t0} /. {y -> y0}]])
MCsample[f_, t_, t0_, y0_, c__, h_] := (Module[{A, tau}, tau = RandomVariate[ExponentialDistribution[1]]; 
   If[tau > t - t0, Return [h*codetofunction[f, c, t0, y0]/E^(-(t - t0))], 
    If[c == {}, Return[MCsample[f, t - tau, t0, y0, {0, 0}, h/E^(-tau)]], 
     If[RandomVariate[UniformDistribution[1]][[1]] <= 0.5, 
       Return[MCsample[f, t - tau, t0, y0, {c[[1]] + 1, c[[2]]}, 2*h/E^(-tau)]],
       A = MCsample[f, t - tau, t0, y0, {0, 0}, 1];
       Return[MCsample[f, t - tau, t0, y0, {c[[1]], c[[2]] + 1}, 2*A*h/E^(-tau)]]]]]])
Solution[f_, t_, t0_, y0_, n_] := (temp = 0; For[i = 1, i <= n, i++, temp += MCsample[f, t, t0, y0, {}, 1]]; Return[temp/n]) 
f[t_, y_] := y^2;Solution[f, 0.5, 0, 1, 100000]
\end{lstlisting}

\noindent
The following Python code implements the algorithm of Theorem~\ref{t1} for systems of ODEs.

\smallskip 

\begin{lstlisting}[language=Python]
import time
import math
import torch
from scipy import special
from torch.distributions.exponential import Exponential
import matplotlib.pyplot as plt
import numpy as np

torch.manual_seed(0)  # set seed for reproducibility

class ODEBranch(torch.nn.Module):
    def __init__(
        self,
        fun,
        t_lo=0.0,
        t_hi=1.0,
        y0=1.0,
        nb_path_per_state=1000000,
        nb_states=6,
        outlier_percentile=1,
        outlier_multiplier=1000,
        patch=1,
        epochs=3000,
        device="cpu",
        verbose=False,
        **kwargs,
    ):
        super(ODEBranch, self).__init__()
        self.fun = fun

        self.loss = torch.nn.MSELoss()
        self.nb_path_per_state = nb_path_per_state
        self.nb_states = nb_states
        self.outlier_percentile = outlier_percentile
        self.outlier_multiplier = outlier_multiplier
        self.patch = patch
        self.t_lo = t_lo
        self.t_hi = t_hi
        self.y0 = y0
        self.dim = len(y0)
        self.epochs = epochs
        self.device = device
        self.verbose = verbose

    def forward(self, code=None):
        start = time.time()
        code = [-1] * self.dim if code is None else code  # start from identity code if not specified
        t = torch.linspace(self.t_lo, self.t_hi, steps=self.nb_states, device=self.device)
        t = t.repeat(self.nb_path_per_state).reshape(self.nb_path_per_state, -1).T
        nb_states_per_patch = math.ceil( ( self.nb_states - 1 ) / self.patch)
        cur_start_idx, cur_end_idx = 0, nb_states_per_patch
        mc_mean, mc_var = [], []
        y0, t0 = torch.tensor(self.y0, device=self.device), torch.tensor(self.t_lo, device=self.device)
        while cur_start_idx < cur_end_idx:
            self.code_to_fun_dict = {}
            t_this_patch = t[cur_start_idx:cur_end_idx]
            H_tensor = torch.ones_like(t_this_patch)
            mask_tensor = torch.ones_like(t_this_patch)
            mc_mean_this_patch = []
            mc_var_this_patch = []
            for i in range(self.dim):
                y = self.gen_sample_batch(
                    t_this_patch,
                    t0,
                    y0,
                    np.array(code),
                    H_tensor,
                    mask_tensor,
                    coordinate=i
                )
                # widen (outlier_percentile, 1 - outlier_percentile) by outlier_multiplier times
                # everything outside this range is considered outlier
                lo = y.nanquantile(self.outlier_percentile/100, dim=1, keepdim=True)
                hi = y.nanquantile(1 - self.outlier_percentile/100, dim=1, keepdim=True)
                lo, hi = lo - self.outlier_multiplier * (hi - lo), hi + self.outlier_multiplier * (hi - lo)
                mask = torch.logical_and(lo <= y, y <= hi)
                mc_mean_this_patch.append((y * mask).sum(dim=1) / mask.sum(dim=1))
                y = y - mc_mean_this_patch[-1].unsqueeze(dim=-1)
                mc_var_this_patch.append(torch.square(y * mask).sum(dim=1) / mask.sum(dim=1))

            # update y0, t0, idx
            mc_mean.append(torch.stack(mc_mean_this_patch))
            mc_var.append(torch.stack(mc_var_this_patch))
            y0, t0 = mc_mean[-1][:, -1], t_this_patch[-1][-1]
            cur_start_idx, cur_end_idx = cur_end_idx, min(cur_end_idx + nb_states_per_patch, self.nb_states)

        if self.verbose:
            print(f"Time taken for the simulations: {time.time() - start:.2f} seconds.")
        return t[:, 0], torch.cat(mc_mean, dim=-1), torch.cat(mc_var, dim=-1)

    @staticmethod
    def nth_derivatives(order, y, x):
        """
        calculate the derivatives of y wrt x with order `order`
        """
        for cur_dim, cur_order in enumerate(order):
            for _ in range(int(cur_order)):
                try:
                    grads = torch.autograd.grad(y.sum(), x, create_graph=True)[0]
                except RuntimeError as e:
                    # when very high order derivatives are taken for polynomial function
                    # it has 0 gradient but torch has difficulty knowing that
                    # hence we handle such error separately
                    # logging.debug(e)
                    return torch.zeros_like(y)

                # update y
                y = grads[cur_dim]
        return y

    def code_to_function(self, code, t, y0, coordinate):
        code = tuple(code)
        if (code, coordinate) not in self.code_to_fun_dict.keys():
            # code (-1, -1, ..., -1) -> identity mapping
            if code == (-1,) * self.dim:
                self.code_to_fun_dict[(code, coordinate)] = y0[coordinate]
            else:
                y = y0.clone().requires_grad_(True)
                self.code_to_fun_dict[(code, coordinate)] = (
                    self.nth_derivatives(code, self.fun(y, coordinate), y).detach()
                )
        return self.code_to_fun_dict[(code, coordinate)]

    def gen_sample_batch(self, t, t0, y0, code, H, mask, coordinate):
        nb_states, _ = t.shape
        tau = Exponential(
            torch.ones(nb_states, self.nb_path_per_state, device=self.device)
        ).sample()
        ans = torch.zeros_like(t)

        ############################### for t + tau >= T
        mask_now = mask.bool() * (t0 + tau >= t)
        if mask_now.any():
            ans[mask_now] = (
                    H[mask_now]
                    * self.code_to_function(code, t0, y0, coordinate)
                    / torch.exp(-(t - t0)[mask_now])
            )

        ############################### for t + tau < T
        mask_now = mask.bool() * (t0 + tau < t)
        if (code == [-1] * self.dim).all():
            if mask_now.any():
                # code (-1, -1,..., -1) -> (0, 0,..., 0)
                tmp = self.gen_sample_batch(
                    t - tau, t0, y0, code + 1, H / torch.exp(-tau), mask_now, coordinate,
                )
                ans = ans.where(~mask_now, tmp)

        else:
            unif = torch.rand(nb_states, self.nb_path_per_state, device=self.device)
            idx = (unif * self.dim).long()
            for i in range(self.dim):
                mask_tmp = mask_now * (idx == i)
                if mask_tmp.any():
                    A = self.gen_sample_batch(
                        t - tau,
                        t0,
                        y0,
                        np.array([0] * self.dim),
                        torch.ones_like(t),
                        mask_tmp,
                        i,
                    )
                    code[i] += 1
                    tmp = self.gen_sample_batch(
                        t - tau,
                        t0,
                        y0,
                        code,
                        self.dim * A * H / torch.exp(-tau),
                        mask_tmp,
                        coordinate,
                    )
                    code[i] -= 1
                    ans = ans.where(~mask_tmp, tmp)
        return ans


if __name__ == "__main__":
    # problem configuration
    device = torch.device("cuda" if torch.cuda.is_available() else "cpu")
    problem = [
        "quadratic",
        "cosine",
        "example_3",
        "example_5",
        "example_6"
    ][4]
    dim = 1
    if problem == "quadratic":
        exact_fun = (lambda t, y, coordinate: y[coordinate] / (1 - y[coordinate] * t))
        f_fun = (lambda y, coordinate: y[coordinate] ** 2)
        t_lo, t_hi = 0, 0.5
        y0 = [1.0] * dim
        nb_states = 6
    elif problem == "cosine":
        def exact_fun(t, y, coordinate):
            return 2 * torch.atan(torch.tanh((t + 2 * math.atanh(math.tan(y[coordinate] / 2))) / 2))
        f_fun = (lambda y, coordinate: torch.cos(y[coordinate]))
        t_lo, t_hi = 0, 1.0
        y0 = [1.0] * dim
        nb_states = 6
    elif problem == "example_3":
        def exact_fun(t, y, coordinate):
            if coordinate == 0:
                return y[coordinate] + t
            else:
                return t + torch.sqrt(y[coordinate] + 2 * t ** 2)

        def f_fun(y, coordinate):
            if coordinate == 0:
                return torch.ones_like(y[0])
            else:
                return (y[coordinate] + y[0]) / (y[coordinate] - y[0])
        t_lo, t_hi = 0, 0.5
        y0 = [t_lo] + [1.0] * dim
        nb_states = 11
    elif problem == "example_5":
        def exact_fun(t, y, coordinate):
            if coordinate == 0:
                return y[coordinate] + t
            else:
                tensor_erfi = (lambda x: special.erfi(x.cpu()).to(device))
                return torch.exp(t**2/2) / (1/y[coordinate] - (math.pi / 2) ** 0.5 * tensor_erfi(t / 2 ** 0.5))

        def f_fun(y, coordinate):
            if coordinate == 0:
                return torch.ones_like(y[0])
            else:
                return y[0] * y[coordinate] + y[coordinate]**2
        t_lo, t_hi = 0, 1.0
        y0 = [t_lo] + [.5] * dim
        nb_states = 11
    elif problem == "example_6":
        def exact_fun(t, y, coordinate):
            if coordinate == 0:
                return t*torch.sin(torch.log(t))
            else:
                return t*torch.cos(torch.log(t))

        def f_fun(y, coordinate):
            if coordinate == 0:
                return ( y[1] + y[0] ) / torch.sqrt(y[0]**2+y[1]**2)
            else:
                return ( y[1] - y[0] ) / torch.sqrt(y[0]**2+y[1]**2)
        t_lo, t_hi = 1, 4
        y0 = [0.0] + [1.0]
        nb_states = 7

    # initialize model and calculate mc samples
    model = ODEBranch(
        f_fun,
        t_lo=t_lo,
        t_hi=t_hi,
        y0=y0,
        device=device,
        nb_states=nb_states,
        verbose=True,
        patch=6,
        outlier_percentile=0.1,
        outlier_multiplier=100,
    )
    t, mc_mean, mc_var = model()
    t_fine = torch.linspace(t_lo, t_hi, 100, device=device)  # finer grid for plotting exact solution
    torch.set_printoptions(precision=5, sci_mode=True)

    # plot exact vs numerical
    for i in range(model.dim):
        print(f"For dimension {i + 1}:")
        print(f"The variance of MC is {mc_var[i]}.")
        print(f"The error squared is {(mc_mean[i] - exact_fun(t, y0, i)) ** 2}.")
        plt.plot(t.cpu(), mc_mean[i].cpu(), '+', label="Numerical solution")
        plt.plot(t_fine.cpu(), exact_fun(t_fine, y0, i).cpu(), label="Exact solution")
        plt.title(f"Dimension {i + 1}")
        plt.legend()
        plt.show()
\end{lstlisting}

\vspace{-0.8cm}

\subsubsection*{Acknowledgement}
 We thank Jiang Yu Nguwi for producing the Python code dealing with systems of ODEs. 

\footnotesize

\newcommand{\etalchar}[1]{$^{#1}$}
\def\cprime{$'$} \def\polhk#1{\setbox0=\hbox{#1}{\ooalign{\hidewidth
  \lower1.5ex\hbox{`}\hidewidth\crcr\unhbox0}}}
  \def\polhk#1{\setbox0=\hbox{#1}{\ooalign{\hidewidth
  \lower1.5ex\hbox{`}\hidewidth\crcr\unhbox0}}} \def\cprime{$'$}

\end{document}